CHAPTER

# Topological Transformation Groups

Alejandro Adem and James F. Davis

*Department of Mathematics, University of Wisconsin, Madison, WI 53706*
*Department of Mathematics, Indiana University, Bloomington, IN 47405*

*Contents*

June 3, 1997 – Draft version; Typeset by LaTeX

Both authors were partially supported by NSF grants.

# 1. Preliminaries

## 1.1. Preface

We will not deal here with the historical background of transformation groups. It suffices to say that they occupy a central rôle in mathematics due to their fundamental importance and ubiquitous nature. Rather we will go straight to the basic objects and examples in the subject and from there describe their development in modern mathematics, emphasizing connections to other areas of algebraic and geometric topology. Indeed it is fairly clear that the basic structure of the theory of transformation groups is complete, hence it is a "mature area" and its most important frontiers lie in the realm of interactions with other areas of mathematics. We simply propose to describe some of the fundamental examples and techniques which make transformation groups an important topic, with the expectation that the interested reader will consult the listed references for a deeper understanding. We also feel that transformation groups continues to be a testing ground for new techniques in algebraic and geometric topology, as well as a source of accessible problems for mathematical research. We thus list some of the basic conjectures still open in the subject, although the interested researchers will be left to find the accessible problems on their own.

Our presentation is organized as follows: in section 1 we deal with basic notions and examples, with the conviction that examples are the best approach for introducing tranformation groups; in section 2 we describe the cohomological aspects associated to group actions which are most relevant in algebraic topology; finally in section 3 we discuss the more geometric aspects of this area. Lists of problems are provided in sections 2 and 3. Finally we would like to make clear that in this text we present a view of transformation groups which reflects our personal interests, omitting such topics as actions of connected Lie groups, and group actions and low-dimensional topology. In no way do we pretend that this is a comprehensive survey of the subject. Points of view on the contents of such a survey will differ; hopefully our list of references will at least point the reader towards other material that may fail to appear in this brief synopsis.

## 1.2. Basic Definitions

A *topological group* is a group which is a Hausdorff topological space, with continuous group multiplication and inversion. Any group can be given the structure of a topological group by equipping the group with the discrete topology. We shall concern ourselves mostly with discrete groups.

A *left action of a topological group $G$ on a Hausdorff space $X$* is a continuous map

$$G \times X \to X$$

$$(g, x) \mapsto gx$$

so that $(gh)x = g(hx)$ and $ex = x$ for all $g, h \in G$ and $x \in X$, where $e \in G$ is

the identity. One says that $X$ is a *G-space*. A *G-map* (or *equivariant map*) is a map $f : X \rightarrow Y$ between $G$-spaces which commutes with the $G$-action, that is, $f(gx) = gf(x)$.

A group action defines a homomorphism

$$\theta : G \rightarrow \mathrm{Homeo}(X)$$

$$g \mapsto (x \mapsto gx)$$

where $\mathrm{Homeo}(X)$ is the group of homeomorphisms of $X$; conversely if $G$ is discrete then any such homomorphism defines a group action. An action is *effective* if $\ker \theta = \{e\}$, that is, for every $g$ there is an $x$ so that $gx \neq x$.

Given a point $x \in X$, define the *orbit* $Gx = \{gx | g \in G\} \subset X$. The *orbit space* $X/G$ is the set of all orbits, given the quotient topology under the obvious surjection $X \rightarrow X/G, \ x \mapsto Gx$. A group action is *transitive* if $X$ consists of a single orbit $Gx$. A typical example of a transitive $G$-space is a homogeneous space $X = G/H$.

Given a point $x \in X$, the *isotropy subgroup* is $G_x = \{g \in G | gx = x\} < G$. Two points in the same orbit have conjugate isotropy groups

$$G_{gx} = gG_x g^{-1}.$$

A group action is *free* if for every point $x \in X$, the isotropy group is trivial, that is, $gx \neq x$ for all $x \in X$ and all $g \in G - \{e\}$. A typical example of a free action is the action of the fundamental group $\pi_1(X, x_0)$ of a connected CW complex on its universal cover $\widetilde{X}$.

An action of a locally compact Hausdorff group $G$ on a space $X$ is *proper* (also termed *properly discontinuous* when $G$ is discrete) if for every $x, y \in X$, there are neighborhoods $U$ of $x$ and $V$ of $y$ so that $\{g \in G | gU \cap V \neq \phi\}$ has compact closure in $G$. If a discrete group acts freely and properly on $X$, then $X \rightarrow X/G$ is a covering space. Conversely if $Y$ is path-connected and has a universal cover $\widetilde{Y}$ and if $H$ is a normal subgroup of $\pi = \pi_1(Y, y_0)$, then $G = \pi/H$ acts freely and properly via deck transformations on $X = \widetilde{Y}/H$ with orbit space $Y$.

### *1.3. Examples*

The subject of transformation groups is motivated by examples. In this section we give various natural examples of group actions on manifolds arising from representation theory and geometry. In later sections we will discuss classification results, regularity results (i.e., to what extent do arbitrary actions resemble naturally occurring ones), and the construction of exotic actions.

By a representation of a topological group $G$, we mean a continuous homomorphism from $G$ to an orthogonal group $O(n)$. Since $O(n)$ acts on a wide variety of spaces, such as $\mathbb{R}^n$, $D^n$, $\mathbb{S}^{n-1}$, $\mathbb{R}P^{n-1}$ and $G_k(\mathbb{R}^n)$, one obtains a multitude of $G$-actions from a representation. Likewise a complex representation $G \rightarrow U(n)$ gives actions on $\mathbb{C}P^{n-1}$, $G_k(\mathbb{C}^n)$, etc. A group action "arising" from a continuous homomorphism $G \rightarrow GL_n(\mathbb{R})$ will be called a *linear action*, however, we won't make that precise. We also remark that any smooth action of a compact Lie group $G$ on a smooth manifold $M$ is locally linear: every $x \in M$ has a neighborhood which is $G_x$-diffeomorphic to a linear $G_x$-action on $\mathbb{R}^n$.

Here are some examples of linear actions. Let $\mathbb{Z}/k = \langle T \rangle$ be a cyclic group of order $k$, let $i_1, \ldots, i_n$ be integers relatively prime to $k$, and let $\zeta_k$ be a primitive $k$-th root of unity. Then $\mathbb{Z}/k$ acts on $\mathbb{S}^{2n-1} \subset \mathbb{C}^n$ via

$$T\Big(z_1, \ldots, z_n\Big) = \Big(\zeta_k^{i_1} z_1, \ldots, \zeta_k^{i_n} z_n\Big).$$

The quotient space $\mathbb{S}^{2n-1}/(\mathbb{Z}/k)$ is the *lens space* $L(k; i_1, \ldots, i_n)$. The quaternion eight group $Q_8 = \{\pm 1, \pm i, \pm j, \pm k\}$ is a subgroup of the multiplicative group of unit quaternions

$$\mathbb{S}^3 = \{a + bi + cj + dk \ \in \ \mathbb{H} \ | a^2 + b^2 + c^2 + d^2 = 1\}$$

and $\mathbb{S}^3/Q_8$ is called the *quaternionic space form*. These are examples of *linear spherical space forms* $\mathbb{S}^{n-1}/G$, which arise from representations $\rho : G \to O(n)$ so that for every $g \in G - \{e\}$, $\rho(g)$ has no $+1$ eigenvalues. The quotient $\mathbb{S}^{n-1}/G$ is then a complete Riemannian manifold with constant sectional curvature equal to $+1$, conversely every such manifold is a linear spherical space form. More generally, complete Riemannian manifolds with constant sectional curvature are called *space forms*, they are quotients of $\mathbb{S}^n$, $\mathbb{R}^n$, or $\mathbb{H}^n$ by a discrete group of isometries acting freely and properly. An excellent discussion is found in Wolf [143].

For a Riemannian manifold $M$ of dimension $n$, the group of isometries is a Lie group, whose dimension is less than or equal to $n(n + 1)/2$; equality is realized only when $M$ is the sphere $\mathbb{S}^n$, real projective space $\mathbb{R}P^n$, Euclidean space $\mathbb{R}^n$, or hyperbolic space $\mathbb{H}^n$. These results are classical, see [80]. If $M$ is compact, so is the isometry group $\mathrm{Isom}(M)$. Conversely if a compact Lie group $G$ acts effectively and smoothly on a compact manifold $M$, then by averaging one can put a Riemannian metric on $M$ so that $G$ acts by isometries. For a closed, smooth manifold $M$, the *degree of symmetry of $M$* is the maximal dimension of a compact Lie group which acts effectively and smoothly on $M$. A systematic study of the degree of symmetry of exotic spheres is found in [70].

Proper actions of infinite discrete groups have been widely studied, especially proper actions on Euclidean spaces. For example, a *crystallographic group* $\Gamma$ is a discrete subgroup of the rigid motions of Euclidean space $\mathrm{Isom}(\mathbb{R}^n)$ so that $\Gamma \backslash \mathrm{Isom}(\mathbb{R}^n)/O(n) = \Gamma \backslash \mathbb{R}^n$ is compact. More generally, a proper action of a discrete group on Euclidean space is determined by a discrete subgroup $\Gamma$ of a Lie group $G$, where $G$ has a finite number of components. Then Iwasawa decomposition theory shows that there is a maximal compact subgroup $K$, unique up to conjugacy, with $G/K$ diffeomorphic to $\mathbb{R}^n$. Given a locally compact group $G$, subgroups $\Gamma$ and $K$ with $\Gamma$ discrete and $K$ compact, then $\Gamma$ acts properly on the homogeneous space $X = G/K$. Suppose $\Gamma$ and $\Gamma'$ are two subgroups of a Lie group $G$, abstractly isomorphic as groups. The question of rigidity [101] asks if they are conjugate subgroups of $G$. The Bieberbach Rigidity Theorem asserts that crystallo-graphic groups are rigid, in the weaker sense that two isomorphic crystallographic groups are conjugate by an affine map of Euclidean space. For many examples of proper actions see [114].

Group actions also play an important part in basic constructions for homotopy theory. Let $X$ denote a topological space with a basepoint: using this point we can obtain natural inclusions $X^n \to X^{n+1}$, where the symmetric group $\Sigma_n$ acts by permutation of coordinates so that these maps are equivariant. The $n$–fold symmetric product on $X$ is defined to be the quotient space $SP^n(X) = X^n/\Sigma_n$, and the infinite symmetric product is defined to be the limit $SP^\infty(X) = \lim_{n \to \infty} SP^n(X)$. A remarkable theorem due to Dold and Thom [52] asserts that $\pi_i(SP^\infty(X)) \cong H_i(X, \mathbb{Z})$. A related construction is the configuration space on $n$ unordered points in $X$, defined by $C_n(X) = (X^n - D)/\Sigma_n$. where $D$ consists of all $n$–tuples $(x_1, \ldots, x_n)$ such that $x_i = x_j$ for some $i \neq j$ (note that the $\Sigma_n$–action is free). These spaces arise in many situations in geometry, topology and physics. In particular if $X = \mathbb{C}$, then $\pi_1(C_n(X)) = B_n$, Artin's braid group on $n$ strings. More sophisticated constructions involving the symmetric groups give rise to models for infinite loop spaces (see [92]).

Covering spaces give natural examples of group actions; we illustrate this with knot theory. If $K$ is a knot (= embedded circle) in $\mathbb{S}^3$, and $n$ is a positive integer, there is a unique epimorphism $\pi_1(\mathbb{S}^3 - K) \to \mathbb{Z}/n$. The corresponding $n$-fold cyclic cover can be completed to a cyclic branched cover $X_n \to \mathbb{S}^3$, that is, $\mathbb{Z}/n$ acts on a closed 3-manifold $X_n$ so that $(X_n/(\mathbb{Z}/n), X_n^{\mathbb{Z}/n}/(\mathbb{Z}/n))$ is $(\mathbb{S}^3, K)$. The homology group $H_1(X_n)$ was the first systematic knot invariant [7], [118] [43].

Exotic (yet naturally occurring) examples of group actions are given by symmetries of Brieskorn varieties [23, Part V, §9]. For a non-zero integer $d$, let $V = V_d^{2n}$ be the complex variety in $\mathbb{C}^{n+1}$ given as the zero set of

$$z_0^d + z_1^2 + \cdots + z_n^2 = 0.$$

The orthogonal group $O(n)$ acts on $V$ fixing the first coordinate and acting on the last $n$ coordinates via matrix multiplication. The variety $V$ has a singularity only at the origin, so

$$\Sigma = \Sigma_d^{2n-1} = V \cap \mathbb{S}^{2n+1}$$

is a smooth $(2n-1)$-dimensional submanifold of $\mathbb{S}^{2n+1}$ and $\Sigma$ is $O(n)$-invariant. Brieskorn investigated the algebraic topology of $\Sigma$ and found that when $n$ and $d$ are both odd, $\Sigma_d^{2n-1}$ is homeomorphic to the sphere, but may have an exotic differential structure. For even $n$, $H_{n-1}(\Sigma) = \mathbb{Z}/d$, and $H_i(\Sigma) = 0$ for $i \neq 0$, $n-1$, $2n-1$. Also $\Sigma_d^3$ is the Lens space $L(d; 1, 1)$. In particular, using the matrix

$$\begin{pmatrix} 1 & 0 & 0 \\ 0 & 1 & 0 \\ 0 & 0 & -1 \end{pmatrix} \in O(3),$$

there is a $\mathbb{Z}/2$-action on $\Sigma_d^5$ with fixed set $\Sigma_d^3 = L(d; 1, 1)$. Since there is no exotic differential structure on $\mathbb{S}^5$, this gives a non-linear $\mathbb{Z}/2$-action on $\Sigma_d^5 \cong \mathbb{S}^5$ for odd $d > 2$. One can also construct non-linear actions on $\mathbb{S}^4$. This stands in contrast to lower dimensions. It is not difficult to show that all actions of finite groups on $\mathbb{S}^1$ and $\mathbb{S}^2$ are homeomorphic to linear actions, and this is conjectured for $\mathbb{S}^3$. It has been shown [99] that all actions of a finite cyclic group on $\mathbb{S}^3$ with fixed set a knot are homeomorphic to a linear action; this was conjectured by P.A. Smith.

### *1.4. Smooth Actions on Manifolds*

A *Lie group* is a topological group which is a smooth ($= C^\infty$) manifold where group multiplication and inversion are smooth maps. A *smooth action of a Lie group $G$ on a smooth manifold $M$* is an action so that $G \times M \to M$ is a smooth map. For a discrete group $G$, there is the corresponding notion of a *PL-action on a PL-manifold*.

The following proposition is clear for discrete $G$, and requires a bit of elementary differential topology [51, II, 5.2] for the general case.

**Proposition 1.1.** *For a smooth, proper, free action of a Lie group $G$ on a manifold $M$, the orbit space $M/G$ admits a smooth structure so that the quotient map $M \to M/G$ is a submersion.*

To make further progress we restrict ourselves to compact Lie groups. To obtain information about $M/G$ we have a theorem of Gleason [60].

**Theorem 1.1.** *Suppose a compact Lie group $G$ acts freely on a completely regular space $X$. Then $X \to X/G$ is a principal $G$-bundle.*

We will give a nice local description (the Slice Theorem) of a smooth action of a compact Lie group. The key results needed are that $G$-invariant submanifolds have $G$-tubular neighborhoods and that orbits $Gx$ are $G$-invariant submanifolds. We will only sketch the theory; for full proofs the reader is referred to Bredon [23] and Kawakubo [78].

**Theorem 1.2.** *Suppose a compact Lie group $G$ acts smoothly on $M$. Any $G$-invariant submanifold $A$ has a $G$-invariant tubular neighborhood.*

**Sketch of Proof:** A $G$-invariant tubular neighborhood is a smooth $G$-vector bundle $\eta$ over $A$ and a smooth $G$-embedding

$$f : E(\eta) \to M$$

onto a open neighborhood of $A$ in $M$ such that the restriction of $f$ to the zero section is the inclusion of $A$ in $M$.

We first claim that $M$ admits a Riemannian metric so that $G$ acts by isometries. By using a partition of unity, one can put an inner product $\ll, \gg$ on the tangent bundle $T(M)$. To obtain a $G$-invariant metric, one averages using the Haar measure on $G$

$$\langle v, w \rangle = \int_G \ll gv, gw \gg dg.$$

Then the exponential map

$$\exp \ : \ W \to M$$

is defined on some open neighborhood $W$ of the zero-section of $T(M)$ by the property that $\exp(X) = \gamma(1)$ where $X \in T_p(M)$ and $\gamma$ is the geodesic so that $\gamma(0) = p$ and $\gamma'(0) = X$. The exponential map is equivariant in the sense that if $X, gX \in W$, then $\exp(gX) = g \exp(X)$. Let $\eta$ be the orthogonal complement of $T(A)$ in $T(M)$, i.e. $\eta$ is the normal bundle of $A$ in $M$. Then one can find a smooth function

$$\epsilon \ : \ A \to \mathbb{R}_{>0}$$

(constant if $A$ is compact) so that

$$\exp \ : \ \overset{\circ}{D}_\epsilon(\eta) \to M$$

is a smooth embedding onto a open neighborhood of $A$ in $M$. The tubular neighborhood is then obtained using a fiber- and zero-section preserving diffeomorphism $E(\eta) \cong \overset{\circ}{D}_\epsilon(\eta)$.

Applying this to the submanifolds $\{x\}$ where $x \in M^G$, one gets

**Corollary 1.1.** *For a smooth action of a compact Lie group $G$ on a manifold $M$, the fixed set $M^G$ is a smooth submanifold.*

Let $x \in M$. The isotropy group $G_x$ is closed in $G$, so is in fact a Lie subgroup. There is a canonical smooth structure on $G/G_x$ so that $\pi : G \to G/G_x$ is a submersion. It is not difficult to show:

**Lemma 1.1.** *Suppose a compact Lie group $G$ acts smoothly on $M$. Let $x \in M$. Then the map $G/G_x \to M$, $g \mapsto gx$ is a smooth embedding. Hence the orbit $Gx$ is a $G$-invariant submanifold of $M$.*

As a corollary of Theorem 1.2 and Lemma 1.1 one obtains:

**Theorem 1.3** (Slice Theorem). *Suppose a compact Lie group $G$ acts smoothly on a manifold $M$. Let $x \in M$. Then there is vector space $V_x$ on which the isotropy group $G_x$ acts linearly and a $G$-embedding*

$$G \times_{G_x} V_x \to M$$

*onto an open set which sends $[g, 0]$ to $gx$.*

For a right $G$-set $A$ and a left $G$-set $B$, $A \times_G B$ is the *Borel construction*, the quotient of $A \times B$ by the diagonal $G$-action. The image of $\{e\} \times V_x$ in $M$ is called a *slice at $x$*. Here the representation $V_x = T_x(Gx)^\perp \subset T_x(M)$, where $G$ acts via isometries of $M$. Then $G \times_{G_x} V$ is diffeomorphic to $T(Gx)^\perp$ and the map in the slice theorem is a $G$-invariant tubular neighborhood of the orbit $Gx$.

We now consider generalizations of the fact that $M^G$ is a smooth submanifold. For a subgroup $H < G$, $M^H$ need not be a manifold. However:

**Theorem 1.4** (Orbit Theorem). *Suppose a compact Lie group acts smoothly on $M$.*

(i) *For a subgroup $H < G$,*

$$M_{(H)} = \{x \in M | H \text{ is conjugate to } G_x\}$$

*is a smooth submanifold of $G$. The quotient map $\pi : M_{(H)} \to M_{(H)}/G$ is a smooth fiber bundle which can be identified with the bundle*

$$G/H \times_{W(H)} (M_{(H)})^H \to (M_{(H)})^H / W(H),$$

*where $W(H) = N(H)/H$ and $N(H)$ is the normalizer of $H$ in $G$.*

(ii) *Suppose $M/G$ is connected. There is an isotropy group $H$ so that for all $x \in M$, $G_x$ is conjugate to a subgroup of $H$. Then $M_{(H)}$ is open and dense in $M$. The quotient $M_{(H)}/G$ is connected.*

Since $G$-invariant submanifolds (e.g. $Gx$, $M^G$, $M_{(H)}$) have $G$-tubular neighborhoods, it behooves us to examine $G$-vector bundles. Recall that a finite-dimensional real representation $E$ of a compact Lie group decomposes into a direct sum of irreducible representations. This decomposition is not canonical, but if one sums all isomorphic irreducible submodules of $E$, then one gets a canonical decomposition. The same thing works on the level of vector bundles.

Let $\mathrm{Irr}(G, \mathbb{R})$ be the set of isomorphism classes of finite-dimensional irreducible $\mathbb{R}G$-modules. For $[V] \in \mathrm{Irr}(G, \mathbb{R})$, let $D(V) = \mathrm{Hom}_{\mathbb{R}G}(V, V)$. Then $D(V)$ equals $\mathbb{R}, \mathbb{C}$, or $\mathbb{H}$.

**Proposition 1.2.** *Let $E$ be a $G$-vector bundle where $G$ is a compact Lie group. Then*

$$\bigoplus_{[V] \in \mathit{Irr}(G, \mathbb{R})} Hom_{\mathbb{R}G}(V, E) \otimes_{D(V)} V \cong E$$

*where the map is $(f, v) \mapsto f(v)$. If $D(V) = \mathbb{C}$ (or $\mathbb{H}$) then the sub-bundle $Hom_{\mathbb{R}G}(V, E) \otimes V$ admits a complex (or symplectic) structure.*

**Corollary 1.2.** *Suppose $\mathbb{Z}/p$ acts smoothly on $M$ with $p$ prime.*

(i) *If $p$ is odd, the normal bundle to the fixed set $M^{\mathbb{Z}/p} \subset M$ admits a complex structure.*

(ii) *If $p = 2$ and the action is orientation-preserving on an orientable manifold $M$, then $\dim M - \dim M^{\mathbb{Z}/2}$ is even.*

For homotopy theoretic information concerning a $G$-space, it is helpful to have the structure of a $G$-CW-complex.

**Definition 1.1.** *A $G$-CW complex is a $G$-space $X$ together with a filtration*

$$\phi = X_{-1} \subset X_0 \subset X_1 \subset X_2 \subset \cdots \subset X_n \subset \cdots \subset X = \bigcup_{n \geqslant 0} X_n$$

*such that $X = colim_{n \to \infty} X_n$ and for any $n \geqslant 0$ there is a pushout diagram*

$$\begin{array}{ccc}
\coprod_{i \in In} G/H_i \times \mathbb{S}^{n-1} & \longrightarrow & X_{n-1} \\
\downarrow & & \downarrow \\
\coprod_{i \in In} G/H_i \times D^n & \longrightarrow & X_n
\end{array}$$

*where $\{H_i\}_{i \in I_n}$ is a collection of subgroups of $G$.*

Another point of view follows. A discrete group $G$ acts cellularly on an ordinary CW complex $X$ if for every $g \in G$ and for every open cell $c$ of $X$, $gc$ is an open cell of $X$ and $gc = c$ implies that $g|_c = \mathrm{Id}$. Any cellular action on a CW complex $X$ gives a $G$-CW-complex and conversely. From this point of view it is clear that if $X$ is a $G$-CW-complex, so are $X/G$ and $X^H$ for all subgroups $H < G$.

Much of the elementary homotopy theory of $CW$-complexes remains valid for $G$-CW-complexes when $G$ is discrete. For example, there is equivariant obstruction theory [24]. Note that specifying a $G$-map $G/H \times \mathbb{S}^{n-1} \to X_{n-1}$ is equivalent to specifying a map $\mathbb{S}^{n-1} \to X_{n-1}^H$. Using this observation, it is easy to show:

**Proposition 1.3** (Whitehead Theorem). *Let $f : X \to Y$ be a $G$-map between $G$-CW-complexes. Then $f$ is a $G$-homotopy equivalence (i.e. there is a $G$-map $g : Y \to X$ so that $f \circ g$ and $g \circ f$ are $G$-homotopic to the identity) if and only if $f^H : X^H \to Y^H$ induces an isomorphism on homotopy groups, for all subgroups $H$ of $G$.*

A smooth $G$-manifold for $G$ a finite group admits an equivariant triangulation, and hence the structure of a $G$-CW-complex [73]. The corresponding result for a smooth, proper action of a Lie group on a manifold appears in [74].

For a smooth $G$-manifold for a finite group $G$, much of the theory of differential topology goes through. For example, there are equivariant Morse functions and equivariant handle decompositions [138]. This leads to equivariant versions of the $s$-cobordism theorem, see [83, Section I.4.C] and the references therein. On the other hand, transversality fails equivariantly: consider the constant $\mathbb{Z}/2$-map $M \longrightarrow \mathbb{R}$ from a manifold with a trivial $\mathbb{Z}/2$-action to the reals with the action $x \mapsto -x$; there is no homotopy to a map which is simultaneously equivariant and transverse.

### *1.4.1. Change of Category*

The subject of actions of groups on *PL* or topological manifolds differs from that of smooth actions on smooth manifolds. An action of a finite group on a topological manifold satisfies none of the regularity theorems of the previous section, and hence have been little studied. For example, one can suspend the involution on $\mathbb{S}^5$ with fixed set $L(d; 1, 1)$ to get an involution on $\mathbb{S}^6$ so that the fixed set (the suspension of the lens space) is not a manifold. Bing [20] constructed an involution on $\mathbb{S}^3$ with fixed set an Alexander horned sphere.

More typically studied are *topologically locally linear actions* of a compact Lie group on a topological manifold or *PL locally linear actions* of a finite group on a *PL* manifold. By definition, these are manifolds with actions which satisfy the conclusion of the Slice Theorem 1.3. Such actions were called locally smooth in the older literature. For such actions the Orbit Theorem 1.4 remains valid; in particular the fixed set $M^G$ is a submanifold. However, equivariant tubular neighborhoods and equivariant handlebodies need not exist. In fact, a locally linear action of a finite group on a closed manifold need not have the *G*-homotopy type of a finite *G-CW*-complex [112]. This makes the equivariant *s*-cobordism theorem [130], [112] in this setting much more subtle; it requires methods from controlled topology. On the other hand, the theory of free actions of finite groups on closed manifolds parallels the smooth theory [79]. For general information on locally linear actions see [23], [140].

### *1.5. Remarks*

In this first section we have introduced basic objects, examples and questions associated to a topological transformation group. In the next section we will apply methods from algebraic topology to the study of group actions. As we shall see, these methods provide plenty of interesting invariants and techniques. After describing the main results obtained from this algebraic perspective, in section 3 we will return to geometric questions. Having dealt with basic cohomological and homotopy–theoretic issues allows one to focus on the essential geometric problems by using methods such as surgery theory. Important examples such as the spherical space form problem will illustrate the success of this approach.

## 2. Cohomological Methods in Transformation Groups

### *2.1. Introduction*

In this section we will outline the important role played by cohomological methods in finite transformation groups. These ideas connect the geometry of group actions to accessible algebraic invariants of finite groups, hence propitiating a fruitful exchange of techniques and concepts, and expanding the relevance of finite transformation groups in other areas of mathematics. After outlining the basic tools in the subject, we will describe the most important results and then provide a selection of topics where these ideas and closely related notions can be applied. Although many results here apply equally well to compact Lie groups, for concreteness we will assume throughout that we are dealing with finite groups, unless stated otherwise. The texts by Allday and Puppe [8], Bredon [23] and tom Dieck [51] are recommended as background references.

To begin we recall a classical result due to Lefschetz: let $X$ be a finite polyhedron and $f : X \to X$ a continuous mapping. The Lefschetz number $L(f)$ is defined as $L(f) = \sum_{i=0}^{dim\, X} (-1)^i Tr\, H_i(f)$, where $H_i(f) : H_i(X;\mathbb{Q}) \to H_i(X;\mathbb{Q})$ is the map induced in rational homology. Lefschetz' fundamental fixed–point theorem asserts that if $L(f) \neq 0$, then $f$ has a fixed point, i.e. an $x \in X$ such that $f(x) = x$. In particular this implies that if $G = \mathbb{Z}/n$ acts on an acyclic finite polyhedron $X$, then $X^G \neq \emptyset$. This result depends on the geometry of $X$ as well as on the simple group–theoretic nature of $\mathbb{Z}/n$. How does this basic result generalize to more complicated groups? In the special case when $G$ is a finite $p$–group ($p$ a prime), P. Smith (see [23]) developed algebraic methods for producing fundamental fixed-point theorems of the type mentioned above. Rather than describe Smith Theory in its original form, we will outline the modern version as introduced by A. Borel in [21], [22].

### 2.2. Universal G-spaces and the Borel Construction

Denote a contractible free $G$–space by $EG$; such an object can be constructed functorially using joins, as was first done by Milnor in [93]. This space is often called a *universal G–space* and has the property that its singular chains are a free resolution of the trivial module over $\mathbb{Z}G$. A cellular model of $EG$ can easily be constructed and from now on we will assume this condition. Now the quotient $BG = EG/G$ is a $K(G, 1)$, hence its cohomology coincides with the group cohomology $H^*(G, \mathbb{Z}) = Ext^*_{\mathbb{Z}G}(\mathbb{Z}, \mathbb{Z})$. The space $BG$ is also called the *classifying space* of $G$ due to the fact that homotopy classes of maps into $BG$ from a compact space $Y$ will classify principal $G$–bundles over $Y$ (a result due to Steenrod [129]). If $X$ is a $G$-space, recall the *Borel Construction* on $X$, defined as

$$X \times_G EG = (X \times EG)/G$$

where $G$ acts diagonally (and freely) on the product $X \times EG$. If $X$ is a point, we simply recover $BG$. If $G$ is any non–trivial finite group, then $EG$ is infinite dimensional; hence if $X$ is a $G$–CW complex, $X \times_G EG$ will be an infinite dimensional CW complex. However, if $G$ acts freely on $X$, then the Borel Construction is homotopy equivalent to the orbit space $X/G$. The cohomology $H^*(X \times_G EG, \mathbb{Z})$ is often called the *equivariant cohomology of the G-space $X$*. In homological terms, this cohomology can be identified with the $G$–hypercohomology of the cellular cochains on $X$ (see [37] or [29] for more on this). Let us assume from now on that $X$ is a finite dimensional $G$–CW complex; although in many instances this condition is unnecessary, it does simplify many arguments without being too restrictive. The key fact associated to the object above is that the projection

$$X \times_G EG \to BG$$

is a fibration with fiber $X$, and hence we have a spectral sequence with

$$E_2^{p,q} = H^p(BG, \mathcal{H}^q(X; \mathbb{A})) \Rightarrow H^{p+q}(X \times_G EG; \mathbb{A}).$$

Note that $G$ may act non–trivially on the cohomology of $X$. We are now in a position to explain the key results from Smith Theory. Let $G = \mathbb{Z}/p$; then the inclusion of the fixed–point set $X^G \to X$ induces a map

$$i_G : X^G \times BG \to X \times_G EG$$

with the following property:

$$i_G^* : H^r(X \times_G EG; \mathbb{F}_p) \to H^r(X^G \times BG; \mathbb{F}_p)$$

is an isomorphism if $r > dim X$. To prove this, we consider the $G$–pair $(X, X^G)$ and the relative Borel Construction $(X, X^G) \times_G EG = (X \times_G EG, X^G \times BG)$. The statement above is equivalent to showing that $H^r((X, X^G) \times_G EG; \mathbb{F}_p) = 0$ for $r$ sufficiently large. However, this follows from the fact that the relative co-chain complex $C^*(X, X^G)$ is $G$–free, and hence the relative equivariant cohomology can be identified with the cohomology of the subcomplex of invariants, which vanishes above the dimension of $X$.

Now if $X$ is mod p homologous to a point, then the spectral sequence collapses and looking at high dimensions we infer that $X^G$ is mod p homologous to a point. Now if $G$ is any finite $p$–group, it will always have a central subgroup of order $p$, hence using induction one can easily show

**Theorem 2.1** (Smith). *If a finite $p$–group $G$ acts on a finite dimensional complex $X$ mod p homologous to a point, then $X^G$ is non–empty and is also mod p homologous to a point.*

In contrast, it is possible to construct fixed–point free actions of $\mathbb{Z}/pq$ (where $p, q$ are distinct primes) on $\mathbb{R}^n$ (see [23]). This indicates that $p$–groups play a distinguished part in the theory of group actions, analogous to the situation in group cohomology or representation theory.

If $G = \mathbb{Z}/p$ acts on $X = \mathbb{S}^n$ with a fixed point, the corresponding spectral sequnce will also collapse. The key observation is that the existence of a fixed-point leads to a cross section for the bundle $X \times_G EG \to BG$, and hence no non-zero differentials can hit the cohomology of the base; as there are only two lines the spectral sequence must collapse. Using induction this yields

**Theorem 2.2** (Smith). *If a finite p–group $G$ acts on a finite dimensional complex $X$ mod $p$ homologous to a sphere with a fixed point, then $X^G$ is also mod $p$ homologous to a sphere.*

Much later, Lowell Jones (see [77]) proved a converse to Smith's theorem for actions on disks which goes as follows.

**Theorem 2.3** (Jones). *Any finite $\mathbb{F}_p$–acyclic complex is the fixed–point set of a $\mathbb{Z}/p$–action and thus of any finite $p$–group on some finite contractible complex.*

The spectral sequence used above can also be applied to prove the following basic result (see [22]).

**Theorem 2.4.** *If $G$, a finite p–group acts on a finite dimensional complex $Y$, then*

$$\sum_{i=0}^{dim\,Y} dim H^i(Y; \mathbb{F}_p) \geqslant \sum_{i=0}^{dim\,Y^G} dim H^i(Y^G; \mathbb{F}_p).$$

Obvious examples such as actions on projective spaces can be analyzed using these techniques; for detailed applications we refer to [23], Chapter VII, which remains unsurpassed as a source of information on this topic. To give a flavour of the results there, we describe an important theorem due to Bredon. Let $P^h(n)$ denote a space such that its mod $p$ cohomology is isomorphic to the ring $\mathbb{F}_p[a]/a^{h+1}$, where $a$ is an element of dimension $n$.

**Theorem 2.5** (Bredon). *Suppose that $p$ is prime and that $G = \mathbb{Z}/p$ acts on a finite dimensional complex $X$ with the mod $p$ cohomology of $P^h(n)$. Then either $X^G = \emptyset$ or it is the disjoint union of components $F_1, \ldots, F_k$ such that $F_i$ is mod $p$ cohomologous to $P^{h_i}(n_i)$, where $h+1 = \sum_{i=1}^k (h_i + 1)$ and $n_i \leqslant n$. The number of components $k$ is at most $p$. For $p$ odd and $h \geqslant 2$, $n$ and the $n_i$ are all even. Moreover, if $n_i = n$ for some $i$, then the restriction $H^n(X; \mathbb{F}_p) \to H^n(F_i; \mathbb{F}_p)$ is an isomorphism.*

## 2.3. Free Group Actions on Spheres

Next we consider applications to the spherical space form problem, namely what finite groups can act freely on a sphere? Let us assume that $G$ does act freely on $\mathbb{S}^n$, then examining the spectral sequence as before we note that it must abut to the cohomology of an $n$–dimensional orbit space, hence the differential

$$d_{n+1} : H^k(G, H^n(\mathbb{S}^n; \mathbb{F}_p)) \to H^{k+n+1}(G; \mathbb{F}_p)$$

must be an isomorphism for $k$ positive, and hence the mod $p$ cohomology of $G$ must be *periodic*. From the Kunneth Formula, it follows that $G \neq (\mathbb{Z}/p)^n$ with $n > 1$;

applying this to all the subgroups in $G$, we deduce that *every* abelian subgroup in $G$ is cyclic, hence obtaining another classical result due to P. Smith.

**Theorem 2.6** (Smith). *If $G$ acts freely on $\mathbb{S}^n$, then every abelian subgroup of $G$ is cyclic.*

A finite group has all abelian subgroups cyclic if and only if its mod $p$ cohomology is periodic for all $p$ (see [37]). Groups which satisfy this condition have been classified and their cohomologies have been computed (see [6]). In this context, a natural question arises: does every periodic group act freely on a sphere? The answer is negative, as a consequence of a result due to Milnor [94]:

**Theorem 2.7** (Milnor). *If $G$ acts freely on $\mathbb{S}^n$, then every element of order 2 in $G$ must be central.*

Hence in particular the dihedral group $D_{2p}$ cannot act freely on any sphere. Note that this result depends on the fact that the sphere is a manifold. However, such restrictions do not matter in the homotopy-theoretic context, as the following result due to Swan [131] shows:

**Theorem 2.8** (Swan). *Let $G$ be a finite group with periodic cohomology; then it acts freely on a finite complex homotopy equivalent to a sphere.*

At this point the serious problem of realizing a geometric action must be addressed; this will be discussed at length in section 3. As a preview we mention the theorem that a group $G$ will act freely on some sphere if and only if every subgroup of order $p^2$ or $2p$ ($p$ a prime) is cyclic; these are precisely the conditions found by Smith and Milnor.

Clearly the methods used for spheres can be adapted to look at general free actions, given some information on the cohomology of the group. The following example illustrates this: let $G$ denote the semidirect product $\mathbb{Z}/p \times_T \mathbb{Z}/p-1$, where the generator of $\mathbb{Z}/p-1$ acts via the generator in the units of $\mathbb{Z}/p$. From this it is not hard to show that $H^*(G; \mathbb{Z}_{(p)}) \cong \mathbb{Z}_{(p)}[u]/pu$, where $u \in H^{2(p-1)}(G; \mathbb{Z}_{(p)})$. Now assume that $G$ acts freely on a connected complex $X$, such that the action is *trivial on homology*. From the spectral sequence associated to this action we can infer the following: the dimension of $X$ must be at least $2(p-1) - 1$. If it were less, then no differential in the spectral sequence could hit the generator from the base; and hence $H^{2(p-1)}(X/G; \mathbb{Z}_{(p)}) \neq 0$, a contradiction.

### 2.4. Actions of Elementary Abelian Groups and the Localization Theorem

Let us now assume that $G = (\mathbb{Z}/p)^r$, an elementary abelian $p$-group. Cohomological methods are extremely effective for studying actions of these groups. Perhaps the most important result is the celebrated "Localization Theorem" due to A. Borel and D. Quillen [110]. To state it we first recall that if $x \in H^1(G; \mathbb{F}_p)$ is non–zero, then its Bockstein $\beta(x)$ is a two dimensional polynomial class. Let $0 \neq e \in H^{2(p^r-1)}(G; \mathbb{F}_p)$ denote the product of all the $\beta(y)$, as $y$ ranges over non-zero elements in $H^1(G; \mathbb{F}_p)$.

**Theorem 2.9** (Borel and Quillen). *Let $G = (\mathbb{Z}/p)^r$ act on a finite dimensional complex $X$. Then, if $S$ is the multiplicative system of powers of $e$, the localized map induced by inclusions*

$$S^{-1} H^*(X \times_G EG; \mathbb{F}_p) \to S^{-1} H^*(X^G \times BG; \mathbb{F}_p)$$

*is an isomorphism.*

This result has substantial applications to the theory of finite transformation groups. Detailed results about fixed–point sets of actions on spheres, projective spaces, varieties, etc. follow from this, where in particular information about the ring structure of the fixed-point set can be provided. An excellent source of information on this is the recent text by Allday and Puppe [8]. An important element to note is that the action of the Steenrod Algebra is an essential additional factor which can be used to understand the fixed–point set (see also [55]). Also one should keep in mind the obvious interplay between the $E_2$ term of the spectral sequence described previously and the information about the $E_\infty$ term the localization theorem provides. Important results which should be mentioned are due to Hsiang [69] and Chang–Skelbred [38]. In particular we have the following fundamental result.

**Theorem 2.10** (Chang and Skelbred). *Let $G = (\mathbb{Z}/p)^r$, then if $X$ is a finite dimensional $G$–CW complex which is also a mod $p$ Poincaré Duality space, then each component $F_i$ of $X^G$ is also a mod $p$ Poincaré Duality space.*

For the case of actions of compact Lie groups, Atiyah and Bott [13] describe a De Rham version of the localization theorem, which is quite useful for studying questions in differential geometry and physics (see also [53]). There are also recent applications of localization techniques to problems in symplectic geometry, for example in [76].

## 2.5. The Structure of Equivariant Cohomology

We now turn to describing qualitative aspects of equivariant cohomology which follow from isotropy and fixed point data. This was originally motivated by attempts to understand the asymptotic growth rate (Krull Dimension) of the mod $p$ cohomology of a finite group $G$. Atiyah and Swan conjectured that it should be precisely the rank of $G$ at $p$ (i.e. the dimension of its largest $p$–elementary abelian subgroup). This result was in fact proved by D. Quillen [110] in his landmark work on cohomology of groups. First we need some notation. Denote by $\mathcal{A}_G$ the family of all elementary abelian $p$–subgroups in $G$, and by $\mathcal{A}_G(X)$ the ring of families $\{f_A : X^A \to H^*(A; \mathbb{F}_p)\}_{A \in \mathcal{A}_G}$ of locally constant functions compatible with respect to inclusion and conjugation. Consider the homomorphism $H^*(X \times_G EG; \mathbb{F}_p) \to \mathcal{A}_G(X)$ which associates to a class $u$ the family $(\tilde{u}_A)$, where $(\tilde{u}_A)$ is the locally constant function whose value at $x$ is the image of $u$ under the map in equivariant cohomology associated to the inclusion $A \subset G$ and the map from a point to $X$ with image $\{x\}$.

**Theorem 2.11** (Quillen). *If $X$ is compact, then the homomorphism above is an $F$–isomorphism of rings, i.e. its kernel and cokernel are both nilpotent.*

The following two results follow from Quillen's work

**Proposition 2.1.** *Let $G$ act on a finite complex $X$ and denote by $p(t)$ the Poincaré series for the mod $p$ equivariant cohomology of $X$. Then $p(t)$ is a rational function of the form $z(t)/\prod_{i=1}^{n}(1 - t^{2i})$, where $z(t) \in \mathbb{Z}[t]$, and the order of the pole of $p(t)$ at $t = 1$ is equal to the maximal rank of an isotropy subgroup of $G$.*

**Proposition 2.2.** *If $G$ is a finite group, then the map induced by restrictions*

$$H^*(G; \mathbb{F}_p) \to lim_{A \in \mathcal{A}_G} H^*(A; \mathbb{F}_p)$$

*is an $F$–isomorphism.*

For example, if $G = \mathcal{S}_n$, the finite symmetric group, then the map above is actually an isomorphism for $p = 2$. We refer the reader to the original paper for complete details; it suffices to say that the proof requires a careful consideration of the Leray spectral sequence associated to the projection $X \times_G EG \to X/G$.

This result has many interesting consequences; here we shall mention that it was the starting point to the extensive current knowledge we have in the cohomology of finite groups (see [6]). An analogous theorem for modules has led to the theory of complexity and many connections with modular representations have been uncovered (see [33] and [17]).

**Example 2.1.** The following simple example ties in many of the results we have discussed. Let $G = Q_8$, the quaternion group of order 8. Its mod 2 cohomology is given by (see [37], [6])

$$H^*(G; \mathbb{F}_2) \cong \mathbb{F}_2[x_1, y_1, u_4]/x_1^2 + x_1 y_1 + y_1^2, x_1^2 y_1 + x_1 y_1^2.$$

Note that the asymptotic growth rate of this cohomology is precisely one, which corresponds to the fact that it is periodic. In addition every element of order 2 is central; in fact $Q_8 \subset \mathbb{S}^3$ and hence acts freely on it by translation. The class $u_4$ is polynomial, transgressing from the top dimensional class in $\mathbb{S}^3$. In fact one can see that

$$H^*(G; \mathbb{F}_2)/(u_4) \cong H^*(\mathbb{S}^3/Q_8; \mathbb{F}_2)$$

which means that the classes $1, x_1, y_1, x_1 y_1, x_1^2 y_1$ represent a cohomology basis for the mod 2 cohomology of the 3–manifold $\mathbb{S}^3/Q_8$. The unique elementary abelian subgroup is the central $\mathbb{Z}/2$, and the four dimensional class $u_4$ restricts to $e_1^4 \in H^4(\mathbb{Z}/2; \mathbb{F}_2)$, where $e_1$ is the one dimensional polynomial generator. The other cohomology generators are nilpotent.

Another interesting group which acts freely on $\mathbb{S}^3$ is the binary icosahedral group $B$ of order 120 (it is a double cover of the alternating group $\mathcal{A}_5$). In this case we have

$$H^*(B; \mathbb{F}_2) \cong \Lambda(x_3) \otimes \mathbb{F}_2[u_4]$$

where as before in the spectral sequence for the group action the top class in the sphere transgresses to $u_4$. From this we obtain $H^*(\mathbb{S}^3/B; \mathbb{F}_2) \cong \Lambda(x_3)$. This orbit space is the Poincaré sphere.

These examples illustrate how geometric information is encoded in the cohomology of a finite group, a notion which has interesting algebraic extensions (see [19]).

### 2.6. Tate Cohomology, Exponents and Group Actions

The cohomology of a finite group can always be computed using a free resolution of the trivial $\mathbb{Z}G$ module $\mathbb{Z}$. It is possible to splice such a resolution with its dual to obtain a *complete resolution* (see [6]), say $\mathcal{F}_*$, indexed over $\mathbb{Z}$, with the following properties: (1) each $\mathcal{F}_i$ is free, (2) $\mathcal{F}_*$ is acyclic and (3) $\mathcal{F}_*, * \geqslant 0$ is a free resolution of $\mathbb{Z}$ in the usual sense.

Now let $X$ be a finite dimensional $G$–CW complex; in [132] Swan introduced the notion of equivariant Tate cohomology, defined as

$$\widehat{H}_G^i(X) \;=\; H^i(Hom_G(\mathcal{F}_*; C^*(X))).$$

An important aspect of the theory is the existence of two spectral sequences abutting to the Tate cohomology, with respective $E_1$ and $E_2$ terms

$$E_1^{p,q} = \widehat{H}^q(G; C^p(X))$$

and

$$E_2^{p,q} = \widehat{H}^p(G; H^q(X; \mathbb{Z}))$$

which arise from filtering the associated double complex in the two obvious ways. Using the fact that free modules are Tate–acyclic, the first spectral sequence can be used to show that $\widehat{H}_G^*(X) \equiv 0$ if and only if the $G$–action is free. More generally one can show that equivariant Tate cohomology depends only on the singular set of the action. In addition it is not hard to see that $\widehat{H}_G^*(X) \cong H^*(X \times_G EG; \mathbb{Z})$ for $* > dim\, X$. Recent work has concentrated on giving a homotopy–theoretic definition of this concept and defining analogues in other theories (see [5], [62]). This involves using a geometric construction of the transfer. Another important ingredient is the 'homotopy fixed point set' defined as $Map_G(EG, X)$; in fact an analysis of the natural map $X^G = Map_G(*, X) \to X^{hG}$ is central to many important results in equivariant stable homotopy.

Let $A$ be a finite abelian group; we define its exponent $\exp(A)$ as the smallest integer $n > 0$ such that $n.a = 0$ for all $a \in A$. Using the transfer, it is elementary to verify that $|G|$ annihilates $\widehat{H}_G^*(X)$; hence exponents play a natural role in this theory. Assume that $X$ is a connected, free $G$–CW complex. Now consider the $E_r^{0,0}$ terms in the second spectral sequence described above; the possible differentials involving it are of the form

$$E_{r+1}^{-r-1,r} \to E_{r+1}^{0,0} \to E_{r+2}^{0,0}$$

with $r = 1, 2, \ldots, N$, $N = \dim X$. From these sequences we obtain that $\exp E_{r+1}^{0,0}$ divides the product of $\exp E_{r+1}^{-r-1,r}$ and $\exp E_{r+2}^{0,0}$ and hence as $E_\infty^{0,0} = 0$, and $E_2^{0,0} = \widehat{H}^0(G; \mathbb{Z}) = \mathbb{Z}/|G|$, we obtain the following condition, first proved by Browder (see [26], [1]): $|G|$ divides the product $\prod_{r=1}^{\dim X} \exp \widehat{H}^{-r-1}(G; H^r(X; \mathbb{Z}))$. We note the following important consequence of this fact.

**Theorem 2.12** (Browder). *If $X$ is a connected, free $(\mathbb{Z}/p)^r$–CW complex, and if the action is trivial in homology, then the total number of dimensions $i > 0$ such that $H^i(X; \mathbb{Z}_{(p)}) \neq 0$ must be at least $r$.*

**Corollary 2.1** (Carlsson). *If $(\mathbb{Z}/p)^r$ acts freely and cellularly on $(\mathbb{S}^n)^k$ with trivial action in homology, then $r \leqslant k$.*

This corollary, was proved by G. Carlsson [34] using different methods. In [4] the hypothesis of homological triviality was removed for odd primes and hence we have the generalization of Smith's result, namely

**Theorem 2.13** (Adem–Browder). *If $p$ is an odd prime and $(\mathbb{Z}/p)^r$ acts freely on $(\mathbb{S}^n)^k$, then $r \leqslant k$.*

For $p = 2$ the same result will hold provided $n \neq 1, 3, 7$. This is a Hopf invariant one restriction.

Another consequence of Browder's result concerns the exponents carried by the Chern classes of a faithful unitary representation of $G$.

**Corollary 2.2.** *Let $\rho : G \to U(n)$ denote a faithful unitary representation of a finite group $G$. Then $|G|$ must divide the product $\prod_{i=1}^n \exp(c_i(\rho))$.*

Using methods from representation theory, one can in fact show [2] that for $G = (\mathbb{Z}/p)^r$,

$$\exp \widehat{H}_G^*(X) = \exp \widehat{H}_G^0(X) = \max\{|G_x|, x \in X\}$$

hence in particular we obtain for any $G$

**Theorem 2.14.** *The Krull Dimension of $H^*(X \times_G EG; \mathbb{F}_p)$ is equal to the maximum value of $log_p\{exp \hat{H}_E^0(X)\}$ as $E$ ranges over all elementary abelian $p$ subgroups of $G$.*

This shows the usefulness of equivariant Tate cohomology, as it will determine asymptotic cohomological information for ordinary equivariant cohomology from a single exponent.

In [27], Browder defined the degree of an action as follows. Let $G$ act on a closed oriented manifold $M^n$ preserving orientation, and let $j : M \to M \times_G EG$ denote the fiber inclusion. Then

$$\deg_G(M) = |H^n(M; \mathbb{Z})/imj^*|.$$

This was independently defined by Gottlieb in [61]; they both show that if $G = (\mathbb{Z}/p)^r$, then $log_p \deg_G(M)$ is equal to the co–rank of the largest isotropy subgroup in $G$. Note in particular that the action will have a fixed–point if and only if $\deg_G(M) = 1$. Using duality it is possible to prove their result from the previous theorem, we refer the reader to [2] for details.

### 2.7. Acyclic Complexes and the Conner Conjecture

If $X$ is a $G$ space and $H \subset G$ is a subgroup then a basic construction is the *transfer map* $C_*(X/G) \to C_*(X/H)$. By averaging on cochains it is elementary to construct such a map (see [29]) with the property that composed with the projection $X/H \to X/G$ the resulting map is multiplication by $[G : H]$ on $H^*(X/G, \mathbb{A})$, where $\mathbb{A}$ is any coefficient group. Note in particular that if $P = \mathrm{Syl}_p(G)$, we have an embedding $H^*(X/G; \mathbb{F}_p) \to H^*(X/P; \mathbb{F}_p)$. A basic result is

**Theorem 2.15.** *If $X$ is a finite dimensional acyclic $G$–complex, then $X/G$ is acyclic.*

From the above, to show that $X/G$ is acyclic it suffices to show (for any prime $p$) that if $X$ is mod $p$ acyclic and $\mathbb{Z}/p$ acts on $X$, then $X/\mathbb{Z}/p$ is mod $p$ acyclic. Consider the the mod $p$ equivariant cohomology of the relative cochain complex for the pair $(X, X^{\mathbb{Z}/p})$; as it is free, we can identify it with the mod $p$ cohomology of the quotient pair, $(X/\mathbb{Z}/p, X^{\mathbb{Z}/p})$. Now the $E_2$ term of the spectral sequence converging to this is of the form $H^p(G, H^q(X, X^{\mathbb{Z}/p}; \mathbb{F}_p))$; using the fact that the fixed–point set must be mod $p$ acyclic (by Smith's theorem) we conclude that it must be identically zero and hence $H^*(X/\mathbb{Z}/p, \mathbb{F}_p) \cong H^*(X^{\mathbb{Z}/p}, \mathbb{F}_p)$ and so $X/\mathbb{Z}/p$ is mod $p$ acyclic. Less obvious is the fact that if $X$ is contractible, then so is $X/G$ (see [51], page 222). The most general results along these lines are due to Oliver [103] who in particular settled a fundamental conjecture due to Conner for compact Lie groups.

**Theorem 2.16** (Oliver)**.** *Any action of a compact Lie group on a Euclidean space has contractible orbit space.*

The main elements in the proof are geometric transfers and a careful analysis of the map $X \times_G EG \to X/G$ which we discussed previously. Oliver also proved some results about fixed–point sets of smooth actions on discs [104], extending a basic example due to Floyd–Richardson (see [23] for details) in a remarkable way.

We introduce a few group–theoretic concepts. Let $\mathcal{G}_p^q$ be the class of finite groups $G$ with normal subgroups $P \triangleleft H \triangleleft G$ such that $P$ is of $p$–power order, $G/H$ is of $q$–power order and $H/P$ is cyclic. Let

$$\mathcal{G}_p = \cup_q \mathcal{G}_p^q, \;\; \mathcal{G} = \cup_p \mathcal{G}_p.$$

We can now state

**Theorem 2.17** (Oliver). *A finite group $G$ has a smooth fixed–point free action on a disk if and only if $G \notin \mathcal{G}$. In particular, any non–solvable group has a smooth fixed–point free action on a disk and an abelian group has such an action if and only if it has three or more non–cyclic Sylow subgroups.*

**Corollary 2.3.** *The smallest abelian group with a smooth fixed–point free action on a disk is $\mathbb{Z}/30 \oplus \mathbb{Z}/30$, of order 900. The smallest group with such an action is the alternating group $\mathcal{A}_5$ of order 60.*

Note that $\mathcal{A}_5$ is precisely the group occurring in the Floyd–Richardson example. Oliver proved a more general version.

**Theorem 2.18.** *For any finite group $G$ not of prime power order, there is an integer $n_G$ (the Oliver number) so that a finite CW-complex $K$ is the fixed–point set of a $G$-action on some finite contractible complex if and only if $\chi(K) \equiv 1$ (mod $n_G$). Furthermore, if $\chi(K) \equiv 1$ (mod $n_G$) there is a smooth $G$-action on a disk with fixed–point set homotopy equivalent to $K$.*

Recently, Oliver [106] has returned to this problem and by analyzing $G$-vector bundles, has determined the possible fixed–point sets of smooth $G$-actions on some disk when $G$ is not a $p$-group.

We now include a small selection of topics in finite transformation groups to illustrate the scope and diversity of the subject, as well as the significance of its applications. This is by no means a complete listing, but hopefully it will provide the reader with interesting examples and ideas.

## 2.8. Subgroup Complexes and Homotopy Approximations to Classifying Spaces

Let $G$ denote a finite group and consider $S_p(G)$, the partially ordered set of all non–trivial $p$–subgroups in $G$. $G$ acts on this object via conjugation and hence on its geometric realization $|S_p(G)|$, which is obtained by associating an n-simplex to a chain of $n + 1$ subgroups under inclusion. Hence we obtain a finite $G$–CW complex inherently associated to any finite group $G$. Similarly if $A_p(G)$ denotes the poset of non–trivial $p$–elementary abelian subgroups, $|A_p(G)|$ will also be a finite $G$–CW complex. These complexes were introduced by K. Brown and then studied by Quillen [111] in his foundational paper. He showed that these complexes have properties analogous to those of Tits Buildings for finite groups of Lie type. Moreover, these geometric objects associated to finite groups are of substantial interest to group theorists, as they seem to encode interesting properties of the group.

We now summarize basic properties of these $G$-spaces.

(1) $|S_p(G)|$ is $G$ equivariantly homotopic to $|A_p(G)|$.
(2) For all $p$–subgroups $P \subset G$, the fixed point set $|S_p(G)|^P$ is contractible.
(3) There is an isomorphism

$$\widehat{H}^*(G; \mathbb{F}_p) \cong \widehat{H}^*_G(|S_p(G)|; \mathbb{F}_p)$$

(due to K. Brown [29]).

(4) In the mod $p$ Leray spectral sequence for the map $|A_p(G)| \times_G EG \to |A_p(G)|/G$ we have that $E_2^{p,q} = 0$ for $p > 0$ and $E_2^{0,q} \cong H^q(G; \mathbb{F}_p)$. This means that $H^*(G; \mathbb{F}_p)$ can be computed from the cohomology of the normalizers of elementary abelian subgroups and their intersections (this is due to P. Webb, see [139] and [6]).

The following example illustrates the usefulness of these poset spaces.

**Example 2.2.** Let $G = M_{11}$, the first Mathieu group. We have that $|A_2(G)|$ is a finite graph, with an action of $G$ on it such that the quotient space is a single edge, with vertex stabilizers $\Sigma_4$ and $GL_2(\mathbb{F}_3)$ and edge stabilizer $D_8$ (dihedral group of order 8). From this information the cohomology of $G$ can be computed (at $p = 2$), and we have (see [6])

$$H^*(G; \mathbb{F}_2) \cong \mathbb{F}_2[v_3, u_4](w_5)/w_5^2 + v_3^2 u_4.$$

Moreover, from the theory of trees we have a surjection

$$\Sigma_4 *_{D_8} GL_2(\mathbb{F}_3) \to G$$

which is in fact a mod 2 cohomology equivalence. Hence the poset space provides an interesting action which in turns leads to a 2–local model for the classifying space of a complicated (sporadic) simple group. More generally this technique can be used to show that if $K$ is a finite group containing $(\mathbb{Z}/p)^2$ but not $(\mathbb{Z}/p)^3$, then at $p$ the classifying space $BK$ can be modelled by using a virtually free group arising from the geometry of the subgroup complex, which is a graph. We refer to [6] for more complicated instances of this phenomenon.

In a parallel development, important recent work in homotopy theory has focused on constructing 'homotopy models' for classifying spaces of compact Lie groups (see [75]). In particular the classifying spaces of centralizers of elementary abelian subgroups can be used to obtain such a model (again $p$–locally). This is related to cohomological results but has a deeper homotopy-theoretic content which we will not discuss here. We suggest the recent paper by Dwyer [54] for a thorough exposition of the homotopy decompositions of classifying spaces. Equivariant methods play an important part in the proofs.

We should also mention that if $G$ is a perfect group, then the homotopy groups $\pi_n(BG^+)$ contain substantial geometric information, often related to group actions. Here $BG^+$ denotes Quillen's *plus construction* which is obtained from $BG$ by attaching two and three dimensional cells and has the property of being simply connected, yet having the same homology as $BG$. We refer the interested reader to [6], Ch. IX for details.

### 2.9. *Group Actions and Discrete Groups*

An important application of finite transformation groups is to the cohomology of discrete groups of finite virtual cohomological dimension, as first suggested by Quillen in [110]. These are groups $\Gamma$ which contain a *finite index* subgroup $\Gamma'$ of finite cohomological dimension (i.e. with a finite dimensional classifying space). Examples will include groups such as amalgamated products of finite groups, arithmetic groups, mapping class groups, etc. If for example $\Gamma \subset GL_n(\mathbb{R})$ is a discrete subgroup, then $\Gamma$ will act on the symmetric space $GL_n(\mathbb{R})/K$ ($K$ a maximal compact subgroup) with *finite* isotropy. Analogous models and their compactifications are the basic building blocks for approaching the cohomology of discrete groups.

More abstractly, using a simple coinduction construction due to Serre (see [29]), one can always build a finite dimensional $\Gamma$–CW complex $X$ such that

(1) $X^H \neq \emptyset$ if and only if $H \subset \Gamma$ is finite,
(2) $X^H$ is contractible for all $H$ finite.

Now we can choose $\Gamma'$ to be a normal subgroup of finite cohomological dimension and finite index in $\Gamma$. Hence the finite group $G = \Gamma/\Gamma'$ will act on the finite dimensional space $X/\Gamma'$, with isotropy subgroups corresponding to the finite subgroups in $\Gamma$. Moreover, it is not hard to see that for a finite subgroup $H \subset \Gamma$,

$$(X/\Gamma')^H \simeq \coprod_{(J)} B(N_\Gamma(J) \cap \Gamma')$$

where $J$ runs over all $\Gamma'$–conjugacy classes of finite subgroups of $\Gamma$ mapping onto $H$ via the projection $\Gamma \to G$ and $N_\Gamma(J)$ is the normalizer of $J$ in $\Gamma$ (see [29]).

We are therefore in an ideal situation to apply Smith Theory to obtain a lower bound on the size of the cohomology of these discrete groups. To make it quite general, we assume given $\Gamma$ of finite cohomological dimension and $P$ a finite $p$–group of automorphisms for $\Gamma$. Let $\overline{\Gamma} = \Gamma \times_T P$, the semi–direct product; now $\Gamma$ is a normal subgroup of finite index in this group. If we choose $J \subset \overline{\Gamma}$ a finite subgroup mapping onto $P$, let $C_\Gamma(J)$ denote its centralizer in $\Gamma$. Let $H^1(P, \Gamma)$ denote the usual non–abelian cohomology and finally denote by $\dim_{\mathbb{F}_p} H^*(Y)$ the total dimension of the homology $\sum H^i(Y; \mathbb{F}_p)$ for a finite dimensional complex $Y$. We can now state (see [3])

**Theorem 2.19.** *If $\Gamma$ is a discrete group of finite cohomological dimension, then for every finite $p$–group of automorphisms $P$ of $\Gamma$ we have*

$$dim_{\mathbb{F}_p} H^*(\Gamma) \geqslant \sum_{J \in H^1(P,\Gamma)} dim_{\mathbb{F}_p} H^*(C_\Gamma(J))$$

*and in particular*

$$dim_{\mathbb{F}_p} H^*(\Gamma) \geqslant dim_{\mathbb{F}_p} H^*(\Gamma^P)$$

*where $\Gamma^P \subset \Gamma$ is the fixed subgroup under the automorphism group $P$.*

As an application of this, we have that if $\Gamma_n(q) \subset SL_n(\mathbb{Z})$ denotes a level $q$ ($q$ prime) congruence subgroup, and if $p$ is another prime, then

$$\dim_{\mathbb{F}_p} H^*(\Gamma(q)) \geqslant 2^{k(p-3)/2} \cdot \dim_{\mathbb{F}_p} H^*(\Gamma_t(q))$$

where $n = k(p-1) + t$, $0 \leqslant t < p - 1$.

The summands in the general formula will represent 'topological special cycles' which in more geometric situation intersect to produce cohomology (see [120]). A result such as the above should be a basic tool for constructing non-trivial cohomology for discrete groups with symmetries; in fact groups such as the congruence subgroups will have many finite automorphisms and hence plentiful cohomology. Equivariant techniques should continue to be quite useful in producing non-trivial cohomology.

We should also mention that K. Brown [29] used equivariant methods to prove very striking results about Euler characteristics of discrete groups. The following is one of them. The group theoretic Euler characteristic of $\Gamma$ (situation as in the beginning of this section) can be defined as $\chi(\Gamma) = \chi(\Gamma')/|G|$; one checks that it is indeed well-defined. Now let $n(\Gamma)$ denote the least common multiple of the orders of all finite subgroups in $\Gamma$. Serre conjectured and K. Brown proved that in fact $n(\Gamma) \cdot \chi(\Gamma) \in \mathbb{Z}$. This beautiful result furnishes information about the size of the finite subgroups in $\Gamma$, provided the Euler characteristic can be computed. In many instances this is the case; for example $\chi(Sp_4(\mathbb{Z})) = -1/1440$, from which we deduce that $Sp_4(\mathbb{Z})$ has subgroups of order 32, 9 and 5. From a more elementary point of view, this result is simply a consequence of the basic fact that the least common multiple of the orders of the isotropy subgroups of a finite dimensional $G$–complex $Y$ (with homology of finite type) must yield an integer when multiplied by $\chi(Y)/|G|$.

## 2.10. Equivariant K–theory

After the usual cohomology of CW complexes was axiomatized by Eilenberg and Steenrod, the introduction of 'extraordinary' theories led to many important results in topology; specifically K–theory was an invaluable tool in solving a number of problems. Atiyah [11] introduced an equivariant version of K–theory whose main properties were developed by Segal [125] and Atiyah-Segal in [14]. We will provide the essential definitions and the main properties which make this a very useful device for studying finite group actions.

Equivariant complex K–theory is a cohomology theory constructed by considering equivariant vector bundles on $G$–spaces. Let $X$ denote a finite $G$–CW complex, a $G$–vector bundle on $X$ is a $G$–space $E$ together with a $G$–map $p : E \to X$ such that

(i) $p : E \to X$ is a complex vector bundle on $X$.

(ii) for any $g \in G$ and $x \in X$, the group action $g : E_x \to E_{gx}$ is a homomorphism of vector spaces.

Assuming that $G$ is a compact Lie group and $X$ is a compact $G$–CW complex then the isomorphism classes of such bundles give rise to an associated Grothendieck group $K_G^0(X)$, which as in the non–equivariant case can be extended to a $\mathbb{Z}/2$ graded theory $K_G^*(X)$, the equivariant complex K–theory of $X$. An analogous theory exists for real vector bundles. We now summarize the basic properties of this theory:

(1) If $X$ and $Y$ are $G$–homotopy equivalent, then $K_G^*(X) \cong K_G^*(Y)$. However, in contrast to ordinary equivariant cohomology, an equivariant map $X \to Y$ inducing a homology equivalence does not necessarily induce an equivalence in equivariant K–theory (see [14]).

(2) $K_G^*(\{x_0\}) \cong R(G)$, the complex representation ring of $G$.

(3) Let $\mathcal{P} \subset R(G)$ denote a prime ideal with support a subgroup $S \subset G$ (in fact $S$ is characterized as minimal among subgroups of $G$ such that $\mathcal{P}$ is the inverse image of a prime of $R(S)$; if $\mathcal{P}$ is the ideal of characters vanishing at $g \in G$, then $S = <g>$), denote by $X^{(S)}$ the set of elements $x \in X$ such that $S$ is conjugate to a subgroup of $G_x$; then we have the following localization theorem due to Segal:

$$K_G^*(X)_{\mathcal{P}} \cong K_G^*(X^{(S)})_{\mathcal{P}}.$$

(4) If $G$ is a finite group, then (see [50])

$$K_G^*(X) \otimes \mathbb{Q} \cong \oplus_{(g)} K^*(X^{<g>}/C_G(g)) \otimes \mathbb{Q}$$

where $g$ varies over all conjugacy classes of elements in $G$. Using this it is possible to identify the Euler characteristic of $K_G^*(X) \otimes \mathbb{Q}$ with the so-called 'orbifold Euler characteristic' [66].

(5) (Completion Theorem, [14])

$$K^*(X \times_G EG) \cong K_G^*(X)\widehat{\phantom{)}}$$

where completion on the right is with respect to the augmentation ideal $I \subset R(G)$ and the module structure arises from the map induced by projection to a point. This is an important result, even for the case when $X$ is a point; it implies that the K–theory of a classifying space can be computed from the completion of the complex representation ring.

We should mention that there is a spectral sequence for equivariant K–theory similar to the Leray spectral sequence discussed before for the projection from the Borel construction onto the orbit space, but which will involve the representation rings of the isotropy subgroups. These basic properties make equivariant K–theory a very useful tool for studying group actions, we refer to [14], [23], [51] for specific applications. The localization theorem ensures that it is particularly effective for actions of cyclic groups. Of course K–theory is also important in index theory [15].

### 2.11. Equivariant Stable Homotopy Theory

Just as in the case of cohomology and K–theory, there is an equivariant version of homotopy theory. In its simplest setting, if $G$ is a finite group and $X$, $Y$ are finite $G$–CW complexes, then we consider $G$–homotopy classes of equivariant maps $f : X \to Y$, denoted $[X, Y]^G$. Such objects and the natural analogues of classical homotopy theoretic results have been studied by Bredon [24] and others, and there is a fairly comprehensive theory. In many instances results are reduced to ordinary homotopy theoretic questions on fixed–point sets, etc. Rather than dwell on this fairly well–understood topic, we will instead describe the basic notions and results in *equivariant stable homotopy theory*, which have had substantial impact in algebraic topology.

Let $V$ denote a finite dimensional real $G$–module and $\mathbb{S}^V$ its 1–point compactification. If $X$ is a finite $G$–CW complex and $Y$ an arbitrary one (both with fixed base points), we can define

$$\{X, Y\}^G \ = \ lim_{U \in \mathcal{U}_G} [\mathbb{S}^U \wedge X, \mathbb{S}^U \wedge Y]^G$$

where $\mathcal{U}_G$ is a countable direct sum of finite dimensional $\mathbb{R}G$–modules so that every irreducible appears infinitely often and the limit is taken over the ordered set of all finite dimensional $G$–subspaces of $\mathcal{U}_G$ under inclusion; and the maps in the directed system are induced by smashing with $\mathbb{S}^{U_1^\perp \cap U_2}$ and identifying $\mathbb{S}^{U_2}$ with $\mathbb{S}^{U_1^\perp \cap U_2} \wedge \mathbb{S}^{U_1}$, where $U_1 \subseteq U_2$. One checks that this is independent of $\mathcal{U}_G$ and identifications using the fact that the limit is attained, by an equivariant suspension theorem.

We can define $\pi_n^G(X) = \{\mathbb{S}^n, X\}^G$ and $\pi_G^n(X) = \{X, \mathbb{S}^n\}^G$, where $X$ is required to be finite in the definition of $\pi_G^n$. The following summarizes the basic properties of these objects.

(1)  (tom Dieck, [51])

$$\pi_*^G(X) \cong \oplus_{(H)} \pi_*^s(EWH_+ \wedge_{WH} X^H)$$

where $WH = NH/H$ and the sum runs over all conjugacy classes of subgroups in $G$.

(2) $\pi_*^G(X)$ and $\pi_G^*(Y)$ are finitely generated for each value of $*$, if $X$ is an arbitrary $G$–complex and $Y$ is a finite $G$-complex.

(3) $\pi_0^G(\mathbb{S}^0) \cong A(G)$ as rings, where $A(G)$ is the Burnside ring of $G$ (see [35]). Note that $\pi_{-*}^G(\mathbb{S}^0) \cong \pi_*^G(\mathbb{S}^0)$ is a module over $A(G)$.

Given the known facts about group cohomology and the complex K–theory of a finite group, it became apparent that the stable cohomotopy of $BG_+$ would be an object of central interest in algebraic topology. Segal conjectured that in dimension zero it should be isomorphic to the I-adic completion of the Burnside ring, an analogue of the completion theorem in K–theory (I the augmentation ideal in $A(G)$). This was eventually proved by Gunnar Carlsson in his landmark 1984 paper (see [35]).

**Theorem 2.20** (Carlsson). *For $G$ a finite group, the natural map $\widehat{\pi_G^*}(\mathbb{S}^0) \to \pi_s^*(BG_+)$ is an isomorphism, where $\widehat{\pi_G^*}(\mathbb{S}^0)$ denotes the completion of $\pi_G^*(\mathbb{S}^0)$ at the augmentation ideal in $A(G)$.*

A key ingredient in the proof is an application of Quillen's work on posets of subgroups to construct a $G$–homotopy equivalent model of the singular set of a $G$–complex $X$ which admits a manageable filtration. The consequences of this theorem have permeated stable homotopy theory over the last decade and in particular provide an effective method for understanding the stable homotopy type (at $p$) for classifying spaces of finite groups (see [90]). For more information we recommend the survey by G. Carlsson [36] on equivariant stable homotopy theory.

This concludes the selected topics we have chosen to include to illustrate the relevance of methods from algebraic topology to finite transformation groups. Next we provide a short list of problems which are relevant to the material discussed in this section.

### 2.12. Miscellaneous Problems

(1) Let $G$ denote a finite group of rank $n$. Show that $G$ acts freely on a finite-dimensional $CW$-complex homotopy equivalent to a product of $n$ spheres $\mathbb{S}^{m_1} \times \ldots \times \mathbb{S}^{m_n}$.

(2) Prove that if $G = (\mathbb{Z}/p)^r$ acts freely on $X = \mathbb{S}^{m_1} \times \ldots \times \mathbb{S}^{m_n}$ then $r \leqslant n$.

(3) Show that if $(\mathbb{Z}/p)^r$ acts freely on a connected CW complex $X$, then

$$\sum_{i=0}^{\dim X} \dim_{\mathbb{F}_p} H^i(X; \mathbb{F}_p) \geqslant 2^r.$$

(4) Find a fixed integer $N$ such that if $G$ is any finite group with $H_i(G; \mathbb{Z}) = 0$ for $i = 1, \ldots, N$, then $G = \{1\}$.

(5) Calculate $K_G^*(|S_p(G)|)$ in representation–theoretic terms.

(6) Show that $|A_p(G)|^G \neq \emptyset$ if and only if there is a normal elementary abelian $p$–subgroup in $G$.

**Remark 2.1.** We have listed only a few, very specific problems which seem directly relevant to a number of questions in transformation groups. Problem (1) would be a generalization of Swan's result, and seems rather difficult. In [18], a solution was provided in the realm of projective $kG$ chain complexes. Problem (2) has been around for a long time and again seems hard to approach. Problem (3) is a conjecture due to G. Carlsson; implies (2) and has an analogue for free chain complexes of finite type. Problem (4) has a direct bearing (via the methods in 2.9) on the problem of (given $G$) determining the minimal dimension of a finite, connected CW complex with a free and homologically trivial action of $G$. Problem (5) is a general formulation of a conjecture due to Alperin in representation theory, as described by Thevenaz [133]. Finally, Problem (6) is a conjecture due to Quillen [111] which has been of some interest in finite group theory (see [9]).

In this section we have attempted to summarize some of the basic techniques and results on the algebraic side of the theory of finite transformation groups. Our emphasis has been to make available the necessary definitions and ideas; additional details can be found in the references. It should however be clear that cohomological methods are a fundamentally useful device for studying transformation groups. In the next section we will consider the more geometric problem of actually constructing group actions when all algebraic restrictions are satisfied; as we will see, the combined approach can be quite effective but unfortunately also rather complicated.

## 3. Geometric Methods in Transformation Groups

The subject of group actions on manifolds is diverse, and the techniques needed for future research seem quite unpredictable, hence we reverse our order of exposition in this section, and start with a discussion of five open problems, the solutions of which would lead to clear advances.

### 3.1. *Five Conjectures*

(i) **Borel Conjecture:** If a discrete group $\Gamma$ acts freely and properly on contractible manifolds $M$ and $N$ with compact quotients, then the quotients are homeomorphic.

(ii) **Group Actions on $\mathbb{S}^3$ are Linear:** Any smooth action of a finite group on $\mathbb{S}^3$ is equivalent to a linear action.

(iii) **Hilbert-Smith Conjecture:** Any locally compact topological group acting effectively on a connected manifold is a Lie group.

(iv) **Actions on Products of Spheres:** If $(\mathbb{Z}/p)^r$ acts freely on $\mathbb{S}^{m_1} \times \cdots \times \mathbb{S}^{m_n}$, then $r \leqslant n$. More generally, what finite groups $G$ act freely on a product of $n$ spheres?

(v) **Asymmetrical manifolds:** There is a closed, simply-connected manifold which does not admit an effective action of a finite group.

### 3.1.1. The Borel Conjecture

It may be a stretch to call the Borel Conjecture a conjecture in transformation groups, but once one has done this, it has to be listed first, as it is one of the main principles of geometric topology. As such, it exerts its influence on transformation groups.

A space is *aspherical* if its universal cover is contractible. The Borel conjecture as stated is equivalent to the conjecture that any two closed, aspherical manifolds with isomorphic fundamental groups are homeomorphic. An aspherical manifold might arise in nature as a complete Riemannian manifold with non-positive sectional curvature or as $\Gamma \backslash G/K$ where $\Gamma$ is a discrete, co-compact, subgroup of a Lie group $G$ with a finite number of components and $K$ is a maximal compact subgroup of $G$, however, the Borel conjecture is a general conjecture about topological manifolds. This is a very strong conjecture; in dimension 3 it implies the Poincaré conjecture, since if $\Sigma^3$ is a homotopy 3-sphere, the conclusion of the Borel conjecture applied to $T^3 \sharp \Sigma^3$ and $T^3 \sharp \mathbb{S}^3$ implies that $\Sigma^3 \cong \mathbb{S}^3$ by Milnor's prime decomposition of 3-manifolds [96]. Nonetheless, the conjecture has been proven in many cases: where one manifold is the $n$-torus $T^n$, $n \geqslant 4$ [72], [136], [79], [58], or if one of the manifolds has dimension $\geqslant 5$ and admits a Riemannian metric of sectional curvature $K \leqslant 0$ [57]. In the study of the Borel conjecture in dimension 3, it is traditional to assume that both manifolds are irreducible, which means that any embedded 2-sphere bounds an embedded 3-ball. This assumption is made to avoid connected sum with a homotopy 3-sphere, and we will call the conjecture that homotopy equivalent, closed, irreducible, aspherical 3-manifolds are homeomorphic the *irreducible Borel conjecture*. The irreducible Borel conjecture has been proven when one of the manifolds is a torus [102], sufficiently large [134], Seifert fibered [124], and work continues in the hyperbolic case [59]. The irreducible Borel conjecture for general hyperbolic 3-manifolds and the Borel conjecture for hyperbolic 4-manifolds remains open.

What is the motivation for the Borel conjecture? First, from homotopy theory – any two aspherical complexes with isomorphic fundamental groups are homotopy equivalent. But the real motivation for Borel's conjecture (made by A. Borel in a coffee room conversation in 1953) was rigidity theory for discrete, co-compact subgroups of Lie groups, in particular the then recent results of Malcev [87] on nilpotent groups and Mostow [100] on solvable groups. Mostow showed that if $\Gamma_1$ and $\Gamma_2$ are discrete, co-compact subgroups of simply-connected solvable Lie groups $G_1$ and $G_2$ (necessarily homeomorphic to Euclidean space), and if $\Gamma_1 \cong \Gamma_2$, then the aspherical manifolds $G_1/\Gamma_1$ and $G_2/\Gamma_2$ are diffeomorphic. In the nilpotent case Malcev showed the stronger statement that there is an isomorphism $G_1 \to G_2$ which restricts to the given isomorphism $\Gamma_1 \to \Gamma_2$. Borel then speculated that while group theoretic rigidity sometimes failed, topological rigidity might always hold. Of course, such phenomena were known prior to the work of Malcev and Mostow. Bieberbach showed rigidity for crystallographic groups. On the other hand, failure of group theoretic rigidity was apparent from the existence of compact Riemann surfaces with the same genus and different conformal structures, i.e there are discrete, co-compact subgroups of $SL_2(\mathbb{R})$ which are abstractly isomorphic, but there is no automorphism of $SL_2(\mathbb{R})$ which carries one to the other. The theory of group theoretic rigidity was investigated further by Mostow [101] and Margulis [88]. The subject of topological rigidity of group theoretic actions (as in Mostow's work on solvable groups) was pursued further by Raymond [117] and his collaborators.

We now discuss variants of the Borel conjecture. The Borel conjecture is not true in the smooth category: smoothing theory shows that $T^n$ and $T^n \sharp \Sigma^n$, $n > 6$ are not diffeomorphic when $\Sigma^n$ is an exotic sphere. The Borel conjecture is not true for open manifolds; there are contractible manifolds not homeomorphic to Euclidean space. This is shown by using the "fundamental group at infinity." In fact, M. Davis [49] constructed closed, aspherical manifolds which are not covered by Euclidean space. There are sharper forms of the Borel conjecture: *a homotopy equivalence between closed, aspherical manifolds is homotopic to a homeomorphism.* There is a reasonable version of the Borel conjecture for manifolds with boundary: *a homotopy equivalence between between compact, aspherical manifolds which is a homeomorphism on the boundary is homotopic, relative to the boundary, to a homeomorphism.*

What should be said for non-free actions? One might call the equivariant Borel conjecture the conjecture that if a discrete group $\Gamma$ acts co-compactly on contractible manifolds $X$ and $Y$ so that the fixed point sets are empty for infinite subgroups of $\Gamma$ and are contractible for finite subgroups of $\Gamma$, then $X$ and $Y$ are $\Gamma$-homeomorphic. This is motivated by the fact that they have the same $\Gamma$-homotopy type. Unfortunately, the equivariant Borel conjecture is not true, however, one can follow the philosophy of S. Weinberger [142] and take the success and failure of the equivariant Borel conjecture in particular cases as a guiding light for deeper investigation.

### 3.1.2. *Group Actions on $\mathbb{S}^3$ are linear*

This is an old question, whose study breaks up into the cases of free and non-free actions. It seems likely that any solution requires geometric input. As is often the case in transformation groups on manifolds, the non-free actions are better understood. In particular, a key case is resolved. P. A. Smith showed that for a prime $p$, if $\mathbb{Z}/p$ acts smoothly, preserving orientation on $\mathbb{S}^3$ with a non-empty fixed point set, then the fixed set is an embedded circle. He conjectured that the fixed set is always unknotted. In [99], it was proven that such an action is equivariantly diffeomorphic to a linear action, giving the Smith conjecture. The proof, building on the work of Thurston, was the joint work of many mathematicians: Bass, Gordon, Litherland, Meeks, Morgan, Shalen, and Yau. The linearization question for general non-free actions is yet unresolved, waiting for a solution for the free case, but linearization results for many non-free actions are given in [99], as well as results concerning the conjecture that any smooth action of a finite group on $\mathbb{R}^3$ is equivalent to a linear action.

The case of free actions is still open, although there has been recent progress. The conjecture may be generalized: a closed 3-manifold with finite fundamental group is diffeomorphic to a linear spherical space form $\mathbb{S}^3/G$. For the trivial group, this is the Poincaré conjecture!

We note that it is not difficult to list the free, linear actions on $\mathbb{S}^3$; the fixed-point free subgroups of $SO(4)$ are given by H. Hopf [68], and then it is easy to give all free representations [143]. Work of Cartan-Eilenberg [37], Milnor [94], Lee [82], Milgram [91], and Madsen [84] gave restrictions on the possible finite fundamental groups of closed 3-manifolds. In another direction, Rubinstein [119] showed that for some groups on Hopf's list (e.g. cyclic groups of small order), any free action on $\mathbb{S}^3$ is equivalent to a linear action. Hamilton [65] showed that a closed 3-manifold with a metric of positive Ricci curvature is diffeomorphic to a linear spherical space form. For recent work on the problem, see [63] where a mixture of gauge theory and surgery theory is used.

Actions of finite groups on $\mathbb{S}^n$, $n \geqslant 4$, are reasonably well understood and need not be equivalent to linear actions. For surveys of non-free actions on $\mathbb{S}^n$, see [99] and [123]. For a survey of free actions on $\mathbb{S}^n$, see [45].

### 3.1.3. Hilbert-Smith Conjecture

We take our discussion of the problem from the surveys of Raymond [116] and Yang [144]. We note at the outset that virtually no progress has been made on this conjecture during the last thirty years, so it may be the time for a fresh look.

The conjecture states that a locally compact topological group $G$ acting effectively on a connected manifold $M$ must be a Lie group. This is known in the following cases:

(i) $G = M$ and the group action is by multiplication. This is the famous result of Montgomery-Zippen which states that a manifold which admits a continuous group structure must be a Lie group.

(ii) $M$ is a differentiable manifold and for all $g \in G$, multiplication by $g$ gives a smooth map $M \to M$. In this case not only is $G$ a Lie group, but the action is also smooth.

(iii) $G$ is compact and every element of $G$ is of finite order. The only such Lie groups are the finite groups.

An inverse limit of finite groups is totally disconnected, hence if the inverse limit is infinite, this gives an example of a compact group which is not a Lie group. The two most obvious examples of such are the infinite $p$-torus $\prod \mathbb{Z}/p$ for a prime $p$ and the additive group $\hat{\mathbb{Z}}_p = \mathrm{invlim}\, \mathbb{Z}/p^n$ of the $p$-adic integers. The infinite $p$-torus cannot act effectively on a manifold by result (iii) above. It is still an open question as to whether the $p$-adic integers can act. In fact, Yamabe has shown that every locally compact group has an open subgroup which is an inverse limit of Lie groups. Using this and (iii) above, one can shown if there is a counterexample to the Hilbert-Smith conjecture, then for some prime $p$, the $p$-adic integers act effectively on a manifold. Such an action would be strange indeed. If $\hat{\mathbb{Z}}_p$ acts effectively on an $n$-dimensional manifold $M$, then $H^{n+2}(M/\hat{\mathbb{Z}}_p; \mathbb{Z}) \neq 0$.

### 3.1.4. Actions on a product of spheres

This problem has been solved when the number of spheres is one; there is the result of Madsen-Thomas-Wall [86] which states that a finite group $G$ acts freely on some sphere if and only if $G$ has no non-cyclic abelian subgroups and no dihedral subgroups. When the number of spheres are greater than one, we discussed algebraic work in section 2, but little geometric work has been done on this problem (but see [105], [67], and [46]). It is evident that *any* finite group will act freely on a product of spheres. Simply take an element $g \in G$, make it act by rotation on an odd sphere and then induce up this action to an action of $G$ on a product of spheres on which $< g >$ still acts freely. Taking products over all elements $g \in G$ provides a product of spheres with a free $G$–action. The main problem which remains unsettled is to show that the number of spheres with any given free $G$–action will bound the rank of the elementary abelian subgroups in the group (see 2.12). On the *constructive* side, the following questions remain unanswered except in some special cases: if $G$ is a finite group of rank $k > 1$ (rank is defined in terms of the maximal $r_p$, taken over all subgroups $(\mathbb{Z}/p)^{r_p}$), does $G$ act freely on a finite dimensional $CW$-complex homotopy equivalent to a product of $k$ spheres? If so, does $G$ act freely on a product of $k$ spheres?

We should mention that certain group–theoretic conditions can be used to produce the required free group actions. For example, if $G$ is a finite 2-group of rank $k$ satisfying Milnor's condition (i.e. every element of order 2 is central) then it will act freely on $(\mathbb{S}^{(|G|/2)-1})^k$. The action is built by inducing up sign representations on $k$ elements of order 2 which span the unique central elementary abelian subgroup in $G$ and then taking their product. More generally it is possible to use this approach to construct actions of arbitrary 2–groups on products of spheres with maximal isotropy of rank equal to the co-rank in $G$ of the largest central elementary abelian subgroup. How to build a free action on a larger product from this object is still unknown. In the context of representation theory the work of U. Ray (see [115]) is also relevant here. She proves that if $G$ is a finite group acting freely on a product of spheres arising from $G$-representations, then the only possible non–abelian composition factors of $G$ are the alternating groups $A_5$ and $A_6$.

### 3.1.5. Asymmetrical manifolds

This problem is not as central as the other problems, but it does point out how little we know about group actions on manifolds not having "obvious" symmetries or manifolds closely resembling such. Presumably the asymmetrical manifold is the generic case (but don't ask what is precisely meant by that!) In the non-simply-connected case, asymmetrical manifolds were first constructed in [41].

### 3.2. Examples and Techniques

#### 3.2.1. Non-linear similarity

A fascinating chapter in the study of transformation groups is topological versus linear similarity. Two linear transformations $T, T' : V \to V$ of a finite dimensional real vector space are *topologically similar* if $T' = hTh^{-1}$ for some homeomorphism $h : V \to V$. Elementary arguments (see [81]) show that if $T$ and $T'$ are topologically similar, then there are decompositions

$$V = V_f \oplus V_\infty \qquad V = V'_f \oplus V'_\infty$$

invariant under $T$ and $T'$ respectively, such that $T \mid_{V_f}$ and $T' \mid_{V'_f}$ have finite (and equal) orders and are topologically similar, while $T \mid_{V_\infty}$ and $T' \mid_{V'_\infty}$ are linearly similar. Thus one may as well assume that $T$ and $T'$ have finite order. It was conjectured that topologically similar implies linearly similar, but this was disproved by Cappell-Shaneson [30] in 1981 using techniques from surgery theory. For $V = \mathbb{R}^9$ and for every $q > 1$, they constructed topologically similar $T$ and $T'$ of order $4q$ which are not linearly similar.

This problem is connected to many others in transformation groups. One first generalizes the problem; two finitely generated $\mathbb{R}G$-modules $V$ and $V'$ for a finite group $G$ are *topologically similar* if there is an equivariant homeomorphism $h : V \to V'$.[1] (Note that if $h$ is differentiable at the origin, then the differential $dh_{\mathbf{0}} : V \to V'$ gives a linear similarity.) The modules $V$ and $V'$ are isomorphic to ones where $G$ acts orthogonally; we assume the actions are orthogonal hereafter. Hence the actions restrict to the unit spheres. The $G$-spaces $S(V)$ and $S(V')$ are *topologically similar* if they are equivariantly homeomorphic. If $S(V)$ and $S(V')$ are topologically similar then so are $V$ and $V'$ (radially extend the homeomorphism) and conversely, if $V$ and $V'$ are topologically similar then $S(V \oplus \mathbb{R})$ and $S(V' \oplus \mathbb{R})$ are topologically similar (one-point compactify). Thus Cappell-Shaneson also constructed examples of non-linearly similar actions on spheres. If the actions on the spheres are free, then Whitehead-Reidemeister-DeRham torsion considerations (see the references in [40]), show that topologically similar actions on $S(V)$ and $S(V')$ are linearly similar. DeRham showed for general linear actions that if the spheres $S(V)$ and $S(V')$ are equivariantly smooth or PL-homeomorphic, then the representations are linearly similar, once again by torsion considerations (see the references in [83]). The fact that non-linear similarities exist implies that equivariant simple homotopy type is not a homeomorphism invariant.

Much of the analysis of the non-linear similarity problem stems from the following observation [31].

---

[1] This notion is connected with foundational issues in study of locally linear actions. If $G$ acts locally linearly on a topological manifold $M$, then every point $x \in M$ has a neighborhood of the form

$$G \times_{G_x} V.$$

The $\mathbb{R}G_x$-module $V$ is only determined up to topological similarity.

**Lemma 3.1.** *$V$ and $V'$ are topologically similar if and only if $S(V)$ and $S(V')$ are G-h-cobordant.*

This means that there is a locally linear $G$-manifold $W$ with boundary $S(V) \coprod S(V')$ so that there are equivariant, orbit-type preserving, strong deformation retracts of $W$ onto $S(V)$ and onto $S(V')$. As part of the definition of $G$-$h$-cobordant, we also require that there is an inverse $G$-$h$-cobordism $-W$ so that $(-W) \cup_{S(V')} W$ and $W \cup_{S(V)} (-W)$ are $G$-homeomorphic rel $\partial$ to $S(V) \times I$ and $S(V') \times I$ respectively.

**Sketch of Proof:** If $h : V \to V'$ is a topological similarity, then by re-scaling one may assume $h(D(V)) \subset \text{int } D(V')$. Let $W = D(V') - \text{int } D(V)$.

Conversely if $W$ is a $G$-$h$-cobordism then

$$\cdots \cup ((-W) \cup W) \cup ((-W) \cup W) \cup \cdots \cong \cdots \cup (-W) \cup (W \cup (-W)) \cup (-W) \cup \cdots$$

Thus $S(V) \times \mathbb{R} \cong S(V') \times \mathbb{R}$. By adding on a point $\{+\infty\}$ to each to compactify one of the ends, we obtain our desired topological similarity.

Construction of non-linearity similarities then proceeds using surgery theory and delicate algebraic number theory. Later approaches use the theory of topological equivariant $h$-cobordisms [31] or bounded methods [64].

The clearest positive result is:

**Theorem 3.1.** *If $G$ has odd order, then topologically similar representations are linearly similar.*

This "Odd Order Theorem" has four different proofs. The result is due independently to to Hsiang-Pardon [71] using stratified pseudo-isotopy theory and lower $K$-groups and to Madsen-Rothenberg [85] using equivariant smoothing theory. Later proofs were given by Rothenberg-Weinberger [121] using Lipschitz analysis and by Hambleton-Pedersen [64] using bounded methods.

Finally we would like to mention a related problem discussed by J. Shaneson in [126]. An $G$-action on a sphere $\Sigma$ is said to be of *Smith type* if for all subgroups $H$, the fixed set $\Sigma^H$ of $H$ is either discrete or connected.

*Given a smooth $G$-action on $\Sigma$ of Smith type so that $\Sigma^G = \{x, y\}$, are the representations of $G$ on the tangent spaces $T_x\Sigma$ and $T_y\Sigma$ linearly similar?*

This question and the similarity question have the same flavor. Given an action on $\Sigma$ as above, then $\Sigma$ minus invariant open disks surrounding $x$ and $y$ gives a candidate for a $G$-$h$-cobordism between $G$-spheres. Conversely, given a $G$-$h$-cobordism $W$ between $G$-spheres, then a candidate action on a sphere $\Sigma$ is the end-point compactification of

$$\cdots \cup (-W) \cup W \cup (-W) \cup W \cup \cdots$$

The study of the non-linear similarity question and the Smith question run parallel; however with the Smith question there are additional considerations of smoothness. Cappell-Shaneson answered Smith's question in the negative, while Sanchez [122] showed that the answer is affirmative for actions of groups of odd order. Cappell-Shaneson modify Smith's question to the conjecture that the tangent space representations are topologically similar. Earlier, Petrie constructed two fixed-point $G$-actions on spheres, with a 2-point fixed set and non-linearly similar tangent space representations, however, which were not of Smith type; we discuss these actions in a later section.

### 3.2.2. Propagation of group actions

Smith theory gives a connection between a group action and homological information at the order of the group. For example, if a group $G$ acts semifreely (i.e. freely away from the fixed set) on a disk or a sphere, then Smith theory shows that the fixed set is a mod $|G|$-disk or sphere. Given a manifold with a group action, propagation is a systematic method for producing group actions on manifolds homologically resembling the given one. Three prototypical questions are:

(i) Given a mod $|G|$-homology sphere $\Sigma$ and a free $G$-action on a sphere, does there exist a free $G$-action on $\Sigma$?

(ii) What are the fixed sets of semifree actions on disks?

(iii) What are the fixed sets of semifree actions on spheres?

Many mathematicians have worked on related ideas; we refer to [140], [47], [32], and [48] for references to original sources. These ideas were pioneered by L. Jones [77] and were taken farthest by S. Weinberger and his collaborators.

For the rest of this section, let $q$ denote the order of the finite group $G$.

**Definition 3.1.** A $G$-action on $Y$ *propagates* across a map $f : X \to Y$ if there exists a $G$-action on $X$ and an equivariant map homotopic to $f$.

Similarly given a $G$-action on $X$ one can talk about propagation across a map $f : X \to Y$.

**Proposition 3.1.** *Let* $f : X \to Y$ *be a map between simply-connected CW-complexes with* $H_*(f; \mathbb{Z}/q) = 0$. *Suppose* $G$ *acts freely on* $Y$, *trivially on* $H_*(Y; \mathbb{Z}[1/q])$. *Then there is a complex* $X'$ *and a homotopy equivalence* $h : X' \to X$ *so that the* $G$-*action propagates across* $f \circ h$. *Furthermore the homotopy type of* $X'/G$ *is uniquely determined.*

We will sketch a proof of the above proposition, to illustrate the homotopy theoretic importance of homological triviality. For a set of primes $P$ and a $CW$-complex $X$, there is a localization map $X \to X_{(P)}$, unique up to homotopy, inducing an isomorphism on $\pi_1$ and a localization of the higher homotopy groups. For an integer $n$, we use the notation $X_{1/n}$ and $X_{(n)}$ to mean invert the primes dividing $n$ and not dividing $n$, respectively. The following lemma is due to Weinberger and is accomplished via a plus construction.

**Lemma 3.2.** *Let $Z$ be a CW-complex with finite fundamental group $G$. Then $G$ acts trivially on $H_*(\widetilde{Z}; \mathbb{Z}[1/q])$ if and only if*

$$Z_{1/q} \simeq \widetilde{Z}_{1/q} \times BG.$$

**Proof of Proposition 3.1:** Let $X'/G$ be the homotopy pullback of

$$
\begin{array}{c}
X_{1/q} \\
\downarrow \\
(Y/G)_{(q)} \longrightarrow X_{(0)} \times BG.
\end{array}
$$

The propagation question can often (in fact, usually) be solved when

$$f : X \to Y$$

is a map between manifolds with $H_*(f; \mathbb{Z}/q) = 0$ and when the $G$-action is trivial on $H_*(\ ; \mathbb{Z}[1/q])$. However the general answer [48], phrased in terms of $K$- and $L$-theory and associated algebraic number theory, is too technical to state here. We give a few examples.

**Example 1:** Let $\Sigma$ be a closed, oriented manifold of dimension $n$, $n \geqslant 5$, $n$ odd, having the $\mathbb{Z}/q$-homology of the sphere. Then any free $G$-action on $\mathbb{S}^n$ can be propagated across any map $f : \Sigma \to \mathbb{S}^n$ whose degree is congruent to 1 modulo $q$.

**Example 2:** Let $\Omega^k \subset D^n$ be a smoothly, properly embedded mod $p$-homology disk ($p$ prime) so that there is a $\mathbb{Z}/p$-action on the normal bundle with fixed set $\Omega^k$. (This happens when the normal bundle admits a complex structure.) Let $N$ be a closed tubular neighborhood. Then the action propagates across the inclusion $\partial N \to D^n - \text{int } N$. Thus there is a smooth, semifree $\mathbb{Z}/p$-action on $D^n$ with fixed set $\Omega^k$.

We next give three theorems giving answers to our three prototypical questions, but note there are many other variations of answers in the literature.

**Theorem 3.2** (Davis-Weinberger [47]). *Let $G$ be a finite group. Let $\Sigma^n$ ($n$ odd, $n \geqslant 5$) be a closed, simply-connected manifold with $H_*(\Sigma^n; \mathbb{Z}/|G|) \cong H_*(\mathbb{S}^n; \mathbb{Z}/|G|)$. Then $G$ acts freely on $\mathbb{S}^n$ if and only if $G$ acts freely on $\Sigma^n$, trivially on homology.*

This produces free actions on manifolds which have no apparent symmetries, as long as they homologically resemble the sphere.

**Theorem 3.3** (Jones, Assadi-Browder, Weinberger [140]). *Let $\Omega^k \subset D^n$ be a proper, smooth embedding with $n - k$ even and greater than 2. Then $\Omega^k$ is the fixed set of a semifree orientation-preserving $G$-action on the disk if and only if*

(i) $\widetilde{H}_*(\Omega^k; \mathbb{Z}/|G|) = 0$

(ii) $\sum (-1)^i [\widetilde{H}_i(\Omega^k; \mathbb{Z})] = 0 \in \widetilde{K}_0(\mathbb{Z}G)$

(iii) *The normal bundle of* $\Omega^k$ *admits a semifree orientation-preserving $G$-action with fixed set* $\Omega^k$.

Here $\widetilde{K}_0(\mathbb{Z}G)$ is the Grothendieck group of finitely generated projective $\mathbb{Z}G$-modules modulo the subgroup generated by free modules. For a finitely generated $\mathbb{Z}G$-module $M$ of finite homological dimension, $[M] \in \widetilde{K}_0(\mathbb{Z}G)$ denotes the Euler characteristic of a projective resolution. In the above theorem, the $\widetilde{K}_0$ obstruction vanishes for $G$ cyclic.

**Theorem 3.4** (Weinberger [141]). *Let $\Sigma^k$ be a PL-locally flat submanifold of $\mathbb{S}^n$ with $n-k$ even and greater than 2. Then $\Sigma^k$ is the fixed set of a semifree orientation-preserving $G$-action on $\mathbb{S}^n$ if and only if $\Sigma^k$ is a $\mathbb{Z}/|G|$-homology sphere, $\mathbb{S}^{n-k-1}$ admits a free linear $G$-action, and certain purely algebraically describable conditions hold for the torsion in the homology of $\Sigma$.*

One condition is condition (ii) from the previous theorem, but there are further conditions involving the Swan subgroup of $L$-theory.

Recently D. Chase [39] has made some progress on propagation in the non-homologically trivial case, but with the presence of addition geometric hypotheses. For example, he has shown the following.

**Theorem 3.5.** *A simply-connected mod 2 homology sphere $\Sigma^{2k}$, $2k \geqslant 6$, has a free involution if and only if there exists an orientation-reversing homeomorphism $L : \Sigma \to \Sigma$ so that $L \circ L$ acts unipotently on homology.*

### 3.2.3. Equivariant surgery

Surgery theory is the primary tool for classification of manifolds and for the study of transformation groups. It is discussed in more detail in Section 3.3. One classically applies surgery theory to transformation groups by doing surgery to the free part of a group action. But in this section we briefly mention a more general type of surgery theory, equivariant surgery. Its development is beset with difficulties arising from a lack of equivariant transversality and embedding theory. The solution to the embedding difficulties is to assume the "gap hypothesis" – that dim $M_{(H)} \geqslant 5$ for all $H \subset G$ and dim $M_{(K)} \geqslant 2$ dim $M_{(H)} + 1$ for all $K \subset H \subset G$, but this assumption is not very appetizing.

Classical surgery theory studies both the uniqueness problem of classifying manifolds up to homeomorphism and the existence problem of determining when a space has the homotopy type of a manifold, however, equivariant surgery theory bifurcates. The equivariant uniqueness question is studied through the Browder-Quinn isovariant theory [28], which is a special case of Weinberger's theory of surgery on stratified spaces [142]. The equivariant surgery to study the existence question was developed by Petrie [109] and was applied and extended by several others, including Dovermann and Dovermann-Schultz. Two of the successes of Petrie's theory were one-fixed point actions on spheres and Smith equivalence of representations and we will discuss these results.

A smooth $G$-action with one fixed point on $\mathbb{S}^n$ gives a fixed point free action on $D^n$ (delete an open equivariant neighborhood of the fixed point), and such an action in turn gives a fixed point free action on $\mathbb{R}^n$ (delete the boundary of the disk). So we first discuss the easier problems of constructing fixed point free actions on Euclidean space and the closed disk. The first examples of fixed point free actions of finite groups on Euclidean space and the closed disk were due to Floyd and Floyd-Richardson respectively, and examples are discussed in [23]. The characterization of groups which can so act was accomplished by Edmonds-Lee [56] and Oliver [104] in the two cases. The first example of a fixed point free action of a finite group on a sphere was due to E. Stein, and the characterization of which groups can act freely without fixed points on some sphere was an application of equivariant surgery due to Petrie [108].

For a finite group $G$, two real representations $V$ and $W$ are said to be *Smith equivalent* if there is a smooth $G$-action on a sphere $\Sigma$ with fixed point set $\{x, y\}$ so that $V \cong T_x\Sigma$ and $W \cong T_y\Sigma$. The first result is along these lines is the following theorem, proven using elliptic differential operators, along with some algebraic number theory.

**Theorem 3.6** (Atiyah-Bott [12], Milnor [97])**.** *If a compact Lie group $G$ acts smoothly and semifreely on a sphere with two fixed points, then the representations at the two fixed points are linearly isomorphic.*

Using equivariant surgery, Petrie and others constructed many examples of smooth actions of finite groups on spheres with two fixed points, but whose representations are not isomorphic (see [109] and the articles in [123]).

### 3.3. *Free actions on spheres*

The techniques we would like to introduce are the theory of the finiteness obstruction and $K_0$, simple homotopy theory and $K_1$, and surgery theory. Rather than introduce these topics abstractly, we would like to introduce these through a concrete geometric situation, the study of topological spherical space forms, manifolds whose universal cover is a sphere.

#### 3.3.1. *Existence: homotopy theoretic techniques*

The existence question is given a finite group $G$, can it act freely on $\mathbb{S}^n$, and the uniqueness question is what is the classification of manifolds with fundamental group $G$ and universal cover $\mathbb{S}^n$. We first discuss the existence question. We will take $n$ to be odd, since by the Lefschetz fixed point theorem, the only group which can act freely on an even-dimensional sphere is the cyclic group of order 2.

If $G$ acts freely on $\mathbb{S}^n$, then by considering the spectral sequence of the fibration $\mathbb{S}^n \to \mathbb{S}^n/G \to BG$, it is easy to show that $H^{n+1}(G, \mathbb{Z}) \cong \mathbb{Z}/|G|$ and that for any additive generator $\alpha \in H^{n+1}(G, \mathbb{Z})$, that $\cup \alpha : \hat{H}^i(G; M) \to \hat{H}^{i+n+1}(G; M)$ is an isomorphism for all $i$ and for all $\mathbb{Z}G$-modules $M$. ($\hat{H}$ is Tate cohomology). Thus $G$ cannot have a subgroup of the form $\mathbb{Z}/p \times \mathbb{Z}/p$. There is a converse.

**Theorem 3.7.** *(Artin-Tate [37, Chapter XII]) Let $G$ be a finite group. The following are equivalent:*

(i) *All abelian subgroups are cyclic.*

(ii) *Every Sylow $p$-subgroup is cyclic or generalized quaternionic.*

(iii) *For some $n$, $H^{n+1}(G, \mathbb{Z}) \cong \mathbb{Z}/|G|$.*

(iv) *For some $n$, there is an element $\alpha \in H^{n+1}(G, \mathbb{Z})$, so that*

$\cup \alpha : \hat{H}^i(G; M) \overset{\sim}{\longrightarrow} \hat{H}^{i+n+1}(G; M)$ *for all $i$ and for all $\mathbb{Z}G$-modules $M$.*

A group satisfying the above conditions is said to be *periodic* and if $H^{n+1}(G, \mathbb{Z}) \cong \mathbb{Z}/|G|$ then $G$ is said to have *period $n+1$*. The periodic groups have been classified and fall into six families (see e.g. [45]). With regard to 3. and 4., for any finite group $G$ and for any $\alpha \in H^{n+1}(G, \mathbb{Z})$, one can show that $\cup \alpha$ is an isomorphism if and only if $\alpha$ is an additive generator of $H^{n+1}(G, \mathbb{Z})$ and $H^{n+1}(G, \mathbb{Z}) \cong \mathbb{Z}/|G|$.

For a group $G$ of period $n+1$, does $G$ act freely on $\mathbb{S}^n$? Not in general, but it does up to homotopy, as we shall see shortly. More precisely, we will show there exists a *Swan complex of dimension $n$*, an $n$-dimensional $CW$-complex $X$ with $\pi_1 X = G$ and $\tilde{X} \simeq \mathbb{S}^n$. It is *polarized* if one fixes the identification of the fundamental group with $G$ and fixes the orientation, i.e. the homotopy class of the homotopy equivalence $\tilde{X} \to \mathbb{S}^n$. By a *periodic projective resolution of period $n+1$*, we mean an exact sequence of $\mathbb{Z}G$-modules

$$0 \to \mathbb{Z} \to P_n \to \cdots \to P_1 \to P_0 \to \mathbb{Z} \to 0,$$

where the $\mathbb{Z}$'s have trivial $G$-actions, and the $P_i$'s are projective over $\mathbb{Z}G$. An example of such is $0 \to \mathbb{Z} \to C_*(\tilde{X}) \to \mathbb{Z} \to 0$, where $X$ is an $n$-dimensional Swan complex. By splicing periodic projective resolutions using the composite $P_0 \to \mathbb{Z} \to P_n$, can form a projective $\mathbb{Z}G$-resolution of $\mathbb{Z}$, or a complete resolution in the sense of Tate cohomology. It follows that $G$ has period $n+1$. In fact, by mapping a projective resolution of $\mathbb{Z}$ to a periodic projective resolution, one defines the $k$-invariant $k \in H^{n+1}(G, \mathbb{Z})$ of the periodic projective resolution. Cup product with $k$ induces periodicity, and hence $k$ is a generator of the cyclic group $H^{n+1}(G, \mathbb{Z})$. The $k$-invariant determines the polarized homotopy type of a Swan complex and the chain homotopy type of a projective periodic resolution. The following theorem is due to Swan [131].

**Theorem 3.8.** *The following are equivalent:*

(i) *$G$ has period $n+1$.*

(ii) *There is a projective periodic resolution $G$ of period $n+1$, where the $P_i$'s are finitely generated.*

(iii) *There is a projective periodic resolution $G$ of period $n+1$, where the $P_i$'s are free.*

(iv) *There is an $n$-dimensional $CW$-complex $X$ with $\pi_1 X = G$ and $\tilde{X} \simeq \mathbb{S}^n$.*

*Furthermore given any generator $\alpha \in H^{n+1}(G, \mathbb{Z})$, one can construct the resolutions in 2. and 3. and the complex in 4. with k-invariant $\alpha$.*

**Discussion of Proof:** The implication (i) $\Rightarrow$ (ii) is the most difficult, although purely homological. We refer the reader to [137]. Here Wall shows that given a generator $\alpha \in H^{n+1}(G, \mathbb{Z})$ and an exact sequence

$$P_n \to P_{n-1} \to \cdots \to P_0 \to \mathbb{Z} \to 0,$$

with the $P_i$'s projective, one can find a periodic projective resolution

$$0 \to \mathbb{Z} \to P'_n \to P'_{n-1} \oplus P_{n-1} \to \cdots \to P_0 \to \mathbb{Z} \to 0$$

with $\alpha$ as k-invariant.

For (ii) $\Rightarrow$ (iii), given a periodic resolution with the $P_i$'s finitely generated projective, add on $Q_0 \xrightarrow{\text{Id}} Q_0$ (in degrees 1 and 0) where $F_0 = P_0 \oplus Q_0$ is free. Continuing inductively, one obtains

$$0 \to \mathbb{Z} \to P'_n \to F_{n-1} \to \cdots \to F_0 \to \mathbb{Z} \to 0$$

with the $F_i$'s finitely generated free. Choose a complement $Q'_n$ so that $P'_n \oplus Q'_n = F_n$ is free. Next we use the Eilenberg swindle.

$$P'_n \oplus (F_n)^\infty = P'_n \oplus (Q'_n \oplus P'_n) \oplus (Q'_n \oplus P'_n) \oplus \cdots \cong (P'_n \oplus Q'_n) \oplus (P'_n \oplus Q'_n) \oplus \cdots = (F_n)^\infty.$$

Thus we can add on $(F_n)^\infty \xrightarrow{\text{Id}} (F_n)^\infty$ (in degrees $n$ and $n-1$) to the above periodic projective resolution to obtain a periodic free resolution.

(iii) $\Rightarrow$ (i) and (iv) $\Rightarrow$ (iii) we have already discussed.

It remains to show ((i), (ii) and (iii)) $\Rightarrow$ (iv). Build a $K(G,1)$ and let $Y$ be its $(n-1)$-skeleton. Then for any generator $\alpha \in H^{n+1}(G, \mathbb{Z})$, we may find a periodic free resolution

$$0 \to \mathbb{Z} \to F_n \xrightarrow{\partial} F_{n-1} \oplus C_{n-1}(\tilde{Y}) \to \cdots \to C_0(\tilde{Y}) \to \mathbb{Z} \to 0$$

with k-invariant $\alpha$. Then let $X = Y \cup (\bigvee \mathbb{S}^{n-1}) \cup (\bigcup e_n)$ where there is a sphere for each $\mathbb{Z}G$-basis element of $F_{n-1}$ and an $n$-cell for each $\mathbb{Z}G$-basis element of $F_n$. The $n$-cells are attached by using the Hurewicz isomorphism $\pi_{n-1}(Y \vee (\bigvee \mathbb{S}^{n-1})) \cong \partial(F_n)$.

The natural question now is whether the $P_i$'s in a periodic projective resolution can be taken to be simultaneously free and finitely generated, or, equivalently, whether the Swan complex can be taken to be finite. The answer to this question is very subtle, and historically was one of the motivations for algebraic $K$-theory.

Let use review the construction of the Swan complex $X$ in the proof of Theorem 3.8. Starting with $Y = K(G,1)^{n-1}$, one can build a periodic projective resolution

$$0 \to \mathbb{Z} \to P_n \xrightarrow{\partial} F_{n-1} \oplus C_{n-1}(\tilde{Y}) \to \cdots \to C_0(\tilde{Y}) \to \mathbb{Z} \to 0$$

where $Y$ is a finite $CW$-complex of dimension $n-1$, $F_{n-1}$ is finitely generated free, and $P_n$ is finitely generated projective. If $P_n$ were stably free (i.e. the direct sum of $P_n$ and a f.g. free module is free), then we don't need the Eilenberg swindle; we could add on the free module in dimensions $n$ and $n-1$ and construct (via the Hurewicz theorem) a finite Swan complex. There is a converse.

**Lemma 3.3.** *Let* $\{Q_n \to \cdots \to Q_0\}$ *and* $\{Q'_n \to \cdots \to Q'_0\}$ *be chain homotopy equivalent chain complexes of projective modules over a ring $R$. Then*

$$Q_n \oplus Q'_{n-1} \oplus Q_{n-2} \oplus \cdots \cong Q'_n \oplus Q_{n-1} \oplus Q'_{n-2} \oplus \cdots.$$

**Corollary 3.1.** *The Swan complex $X$ has the homotopy type of a finite complex if and only if $P_n$ is stably free.*

**Proof of Lemma 3.3** : First prove it when one complex is zero by induction on $n$. Next prove it in general by applying the acyclic case to the algebraic mapping cone of the chain homotopy equivalence.

This result leads to the notion of the *reduced projective class group* $\tilde{K}_0(R)$ *of a ring $R$.* Elements are represented by $[P]$ where $P$ is a finitely generated projective $R$-module. Here $[P] = [Q]$ if and only if $P \oplus R^m \cong Q \oplus R^n$ for some $m$ and $n$. This is a classical notion; if $R$ is a Dedekind domain, then $\tilde{K}_0(R)$ can be identified with the ideal class group of $R$ (see [98]). The above corollary shows that a Swan complex $X$ defines an element $[X] \in \tilde{K}_0(\mathbb{Z}G)$ (represented by $P_n$) which vanishes if and only if $X$ has the homotopy type of a finite complex.

Next we analyze what happens to the finiteness obstruction if the $k$-invariant is changed. Suppose $X'$ and $X$ are Swan complexes. Then there is a map $X' \to X$ of degree $d$, and the $k$-invariants satisfy $k' = dk$, with $(d, |G|) = 1$. Define $P_d = \ker(\epsilon)$, where $\epsilon: \mathbb{Z}G \to \mathbb{Z}/d$ is defined by $\epsilon(\sum n_g g) = \sum n_g$. Then $P_d$ is projective (in fact, $P_d \oplus P_e$ is free if $de \equiv 1 \pmod{|G|}$). The finiteness obstructions satisfy

$$[X'] = [X] + [P_d].$$

Swan [131] defined what is now called the *Swan subgroup* of $\tilde{K}_0(\mathbb{Z}G)$ as $T(G) = \{[P_d] \mid (d, |G|) = 1\}$. For a group $G$ having period $n+1$, one defines the *Swan finiteness obstruction* $\sigma_{n+1}(G) \in \tilde{K}_0(\mathbb{Z}G)/T(G)$ by setting $\sigma_{n+1}(G) = [X]$ for any Swan complex $X$ of dimension $n$ and fundamental group $G$. Then $\sigma_{n+1}(G) = 0$ if and only if there exists a *finite* Swan complex of dimension $n$ and fundamental group $G$

Here are three fascinating examples (see [45] for more details).

(i) (Swan [131]) There is a finite 3-dimensional Swan complex $X$ for the dihedral group of order 6 (so $\sigma_4(D_6) = 0$). By Milnor [94], $X$ does not have the homotopy type of a manifold.

(ii) (Swan [131], Martinet [89]) There is a 3-dimensional Swan complex with fundamental group the quaternion group of order 8 which does not have the homotopy type of a *finite* complex, and hence does not have the homotopy type of a closed manifold. (There is a 3-dimensional quaternionic space form, thus $\sigma_4(Q_8) = 0$ but $T(Q_8) \neq 0$.)

(iii) (Davis [42], Milgram [91]) Let $G = \mathbb{Z}/3 \times_T Q_{16}$ be a semidirect product where the element of order 8 in the quaternion group of order 16 acts non-trivially on $\mathbb{Z}/3$. $G$ has period 4 but there is no *finite* 3-dimensional Swan complex, and hence no closed 3-manifold with fundamental group $G$ (cf. [82]) (and so $\sigma_4(G) \neq 0$). The above is the smallest group satisfying this and is due to Davis; Milgram was the first to find examples with non-zero finiteness obstruction.

The work of Swan was generalized by Wall to the theory of the Wall finiteness obstruction.

**Definition 3.2.** A *CW*-complex $X$ is *finitely dominated* if there is a finite *CW*-complex $X_f$ and maps $i : X \to X_f$ and $r : X_f \to X$ so that $r \circ i \simeq \mathrm{Id}_X$.

This is similar to saying a module is finitely generated projective if and only if it is the retract of a finitely generated free module.

**Lemma 3.4.** *A Swan complex $X$ is finitely dominated.*

**Proof:** By the proof of Theorem 3.8 we may assume $X$ is homotopy equivalent to $Z = Y \cup (\bigvee \mathbb{S}^{n-1}) \cup (\bigcup e^n)$ where $Y$ is a finite $(n-1)$-dimensional *CW*-complex. Furthermore

$$C_*(\tilde{Z}) = F_n^\infty \oplus P_n \to F_n^\infty \oplus C_{n-1}(\tilde{Y}) \to \cdots \to C_0(\tilde{Y}),$$

where $F_n = P_n \oplus Q_n$ is a finitely generated free module and

$$0 \to \mathbb{Z} \to P_n \xrightarrow{\partial} C_{n-1}(\tilde{Y}) \to \cdots \to C_0(\tilde{Y}) \to \mathbb{Z} \to 0$$

is a periodic projective resolution. Then construct $X_f = Y \cup (\bigcup e^n)$ with

$$C_*(X_f) = P_n \oplus Q_n \xrightarrow{\partial \oplus 0} C_{n-1}(\tilde{Y}) \to \cdots \to C_0(\tilde{Y}),$$

using the Hurewicz isomorphism $\pi_{n-1}Y \cong (\partial \oplus 0)(P_n \oplus Q_n)$ to attach the $n$-cells. There are obvious chain maps $i_* : C_*(\tilde{X}) \to C_*(\tilde{X}_f)$ and $r_* : C_*(\tilde{X}_f) \to C_*(\tilde{X})$ so that $r_* \circ i_*$ is chain homotopic to the identity. Use the relative Hurewicz theorem to extend the inclusions $Y \to X_f$ and $Y \to X$ to the desired maps $i : Z \to X_f$ and $r : X_f \to Z$

Wall [135] generalized this.

**Proposition 3.2.** *A connected CW-complex $X$ is finitely dominated if and only if $\pi_1 X$ is finitely presented and $C_*(\tilde{X})$ is chain homotopy equivalent to a complex $P_k \to \cdots \to P_0$ of finitely generated projective $\mathbb{Z}[\pi_1 X]$-modules.*

**Definition 3.3.** For a finitely dominated connected *CW*-complex $X$, the *Wall finiteness obstruction* $[X] \in \tilde{K}_0(\mathbb{Z}\pi_1 X)$ is defined by $[X] = \Sigma(-1)^i[P_i]$ where $C_*(\tilde{X})$ is chain homotopy equivalent to $\{P_k \to \cdots \to P_0\}$, a chain complex of finitely generated $\mathbb{Z}[\pi_1 X]$-modules.

By Lemma 3.3, $[X]$ is well-defined.

**Theorem 3.9.** *(Wall [135]) Let $X$ be a finitely dominated connected CW-complex. Then $[X] = 0$ if and only if $X$ is homotopy equivalent to a finite CW-complex.*

### 3.3.2. Uniqueness : homotopy theoretic techniques

We now turn to the uniqueness question. The general question is the classification up to homeomorphism of manifolds covered by a sphere; we consider only a very special case, the classification of classical lens spaces $L = L(k; i_1, \ldots, i_n)$. For references on this material see Cohen [40] or Milnor [97]. Recall that $L$ has fundamental group $\mathbb{Z}/k$, its universal cover is $\mathbb{S}^{2n-1} \subset \mathbb{C}^n$, and the integers $i_1, \ldots, i_n$ are relatively prime to $k$; they give the rotations in the complex factors of $\mathbb{C}^n$. It is easy to see that the map $[z_1, \ldots, z_n] \mapsto [z_1^{i_1}, \ldots, z_n^{i_n}]$ gives a map $L(k; 1, \ldots, 1) \to L(k; i_1, \ldots, i_n)$ of degree $\prod_{j=1}^n i_j$. The following proposition follows from the earlier discussion of $k$-invariants.

**Proposition 3.3.** *Two lens spaces $L(k; i_1, \ldots, i_n)$ and $L(k; i'_1, \ldots, i'_n)$ are homotopy equivalent if and only if $i_1 \cdots i_n \equiv \pm a^n i'_1 \cdots i'_n \pmod{k}$ for some integer $a$.*

The plus or minus corresponds with the choice of orientation and the factor of $a^n$ corresponds with the identification of the fundamental group of the lens space with $\mathbb{Z}/k$.

There are some obvious diffeomorphisms between lens spaces. For example, if $(i_1, \ldots, i_n)$ considered modulo $k$ is a permutation of $(i'_1, \ldots, i'_n)$, then permuting the complex coordinates gives a diffeomorphism from $L(k; i_1, \ldots, i_n)$ to $L(k; i'_1, \ldots, i'_n)$. Similarly, mapping $[z_1, \ldots, z_j, \ldots z_n] \mapsto [z_1, \ldots, \overline{z}_j, \ldots, z_n]$ gives a diffeomorphism from $L(k; i_1, \ldots, i_j, \ldots, i_n)$ to $L(k; i_1, \ldots, -i_j, \ldots, i_n)$. Finally, if $a$ is relatively prime to $k$, the map $[z_1, \ldots, z_n] \mapsto [\zeta^a z_1, \ldots, \zeta^a z_n]$ gives a diffeomorphism from $L(k; i_1, \ldots, i_n)$ to $L(k; a i_1, \ldots, a i_n)$. Lens spaces exhibit rigidity; all diffeomorphisms are generated by the above three types.

**Theorem 3.10.** *Let $L(k; i_1, \ldots, i_n)$ and $L' = L(k; i'_1, \ldots, i'_n)$ be two lens spaces. The following are equivalent.*

(i) *$L$ and $L'$ are isometric, where they are given the Riemannian metrics coming from being covered by the round sphere.*

(ii) *$L$ and $L'$ are diffeomorphic.*

(iii) *$L$ and $L'$ are homeomorphic.*

(iv) *$L$ and $L'$ have the same simple homotopy type.*

(v) *There are numbers $\epsilon_1, \ldots, \epsilon_n \in \{1, -1\}$ and a number $a$ relatively prime to $k$, so that $(i_1, \ldots, i_n)$ is a permutation of $(a\epsilon_1 i_1, \ldots, a\epsilon_n i_n)$ modulo $k$.*

The implications $(v) \Rightarrow (i)$, $(i) \Rightarrow (ii)$, and $(ii) \Rightarrow (iii)$ are easy, and $(iii) \Rightarrow (iv)$ follows from a theorem of Chapman which states that homeomorphic finite $CW$-complexes have the same simple homotopy type. To proceed farther, one must introduce the torsion and the Whitehead group. For a ring $R$, embed $GL_n(R)$ in $GL_{n+1}(R)$ by $A \mapsto \begin{pmatrix} A & 0 \\ 0 & 1 \end{pmatrix}$. Let $GL(R) = \cup_{n=1}^{\infty} GL_n(R)$. Define

$$K_1(R) = GL(R)/[GL(R), GL(R)]$$

and $\widetilde{K}_1(R) = K_1(R)/(-1)$ where $(-1) \in GL_1(R)$. If $R$ is a commutative ring, the determinant gives a split surjection $\det : K_1(R) \to R^*$, and if $R$ is a field, this is an isomorphism.

A *based $R$-module* is a free $R$-module with a specified basis. A chain complex $C$ over a ring $R$ is *based* if each $C_i$ is based, *finite* if $\oplus_i C_i$ is a finitely generated $R$-module, and *acyclic* if the its homology is zero. We now assume, for the sake of simplicity of exposition, that $R$ has the property that $R^n \cong R^m$ implies that $n = m$ (e.g. a group ring $\mathbb{Z}\pi$ has this property since it maps epimorphically to $\mathbb{Z}$). An isomorphism $f : M \to M'$ between based $R$-modules determines an element $[f] \in \widetilde{K}_1(R)$. Given an chain isomorphism $f : C \to C'$ between finite, based chain complexes, define the *torsion of $f$* by

$$\tau(f) = \prod [f_i : C_i \to C'_i]^{(-1)^i} \in \widetilde{K}_1(R).$$

We next indicate the definition of the torsion $\tau(C)$ of a finite, based, acyclic chain complex $C$. Implicit in the proof of Lemma 3.3 is an algorithm for computing it.

**Proposition 3.4.** *Let $\mathcal{C}$ be the class of finite, acyclic, based chain complexes over $R$. Then there is a unique map $\mathcal{C} \to \widetilde{K}_1(R)$, $C \mapsto \tau(C)$ satisfying the following axioms*

   (i) *If $f : C \to C'$ is a chain isomorphism where $C, C' \in \mathcal{C}$, then $\tau(C') = \tau(f)\tau(C)$.*

   (ii) $\tau(C \oplus C') = \tau(C)\tau(C')$.

   (iii) $\tau(0 \to C_n \xrightarrow{d} C_{n-1} \to 0) = (-1)^{n-1}[d]$.

We use this to define the torsion of a homotopy equivalence. If $f : C \to C'$ is a chain homotopy equivalence between finite, based $R$-chain complexes, define $\tau(f) = \tau(C(f))$ where $C(f)$ is the algebraic mapping cone of $f$. The reader should check that when $f$ is a chain isomorphism, our two definitions of $\tau(f)$ agree. With somewhat more difficulty, the reader can verify that axiom 1 above holds for a chain homotopy equivalence $f : C \to C'$ between acyclic, finite, based complexes.

Let $X$ be a finite, connected $CW$-complex with fundamental group $\pi$. For every cell in $X$, choose a cell in the universal cover $\widetilde{X}$ lying above it, and an orientation for that cell. Then $C(\widetilde{X})$ is a finite, based $\mathbb{Z}\pi$-chain complex. Accounting for the ambiguity in the choice of basis, define the *Whitehead Group* $Wh(\pi) = K_1(\mathbb{Z}\pi)/\{\pm g\}_{g \in \pi}$ where $\{\pm g\}$ refers to the image of a one-by-one matrix. For a homotopy equivalence $h : X \to Y$, where $X$ and $Y$ are finite, connected $CW$-complexes with fundamental group $\pi$, define the *torsion of* $h$, $\tau(h) \in Wh(\pi)$ by $\tau(h) = \tau(\widetilde{h} : C(\widetilde{X}) \to C(\widetilde{Y}))$. A homotopy equivalence $h$ is a *simple homotopy equivalence* if $\tau(h) = 0$. There is a geometric interpretation: $h$ is a simple homotopy equivalence if and only if $X$ can be obtained from $Y$ via a sequence of elementary expansions and collapses; see [40]. A compact manifold has a canonical simple homotopy type; for smooth manifolds this follows from triangulating the manifold [97] and for topological manifolds this is more difficult [79].

Given a homotopy equivalence between lens spaces

$$h : L = L(k; i_1, \ldots, i_n) \longrightarrow L' = L(k; i'_1, \ldots, i'_n)$$

it is possible to compute its torsion in terms of $(i_1, \ldots, i_n)$ and $(i'_1, \ldots, i'_n)$. We will do something technically easier, which still leads to a proof of of Theorem 3.10. Let $T \in \mathbb{Z}/k$ denote the generator used to define the lens spaces and let $\Sigma = 1 + T + \cdots + T^{k-1} \in \mathbb{Z}[\mathbb{Z}/k]$ denote the norm element. There is a decomposition of rings

$$\mathbb{Q}[\mathbb{Z}/k] \cong \mathbb{Q} \times \Lambda,$$

where $\Lambda = \mathbb{Q}[\mathbb{Z}/k]/\Sigma$.

**Definition 3.4.** If $L$ is a finite complex with fundamental group $\mathbb{Z}/k = \langle T \rangle$, so that $T$ acts trivially on $H_*(\widetilde{L}; \mathbb{Q})$, define the *Reidemeister torsion*

$$\Delta(L) = \tau(C(L; \Lambda)) \in K_1(\Lambda)/\{\pm T\}.$$

Finally, we can outline the proof of (iv) implies (v) in Theorem 3.10. Let $h : L \to L'$ be a simple homotopy equivalence of lens spaces. Then

$$\tau(C(L; \Lambda)) = \tau(C(L'; \Lambda)) \in K_1(\Lambda)/\{\pm T\}.$$

We now wish to compute both sides. Implicit in the definition of a lens space is an identification of the fundamental group with $\mathbb{Z}/k$. By perhaps replacing $L' = L(k; i'_1, \ldots, i'_n)$ by the diffeomorphic space $L(k; ai'_1, \ldots, ai'_n)$, we assume that $h$ induces the identity on the fundamental group.

For $\widetilde{L}(k; i_1) = \mathbb{S}^1$ it is easy to see that

$$C(\widetilde{L}(k; i_1)) \cong (\mathbb{Z}[\mathbb{Z}/k] \xrightarrow{\cdot g^{j_1} - 1} \mathbb{Z}[\mathbb{Z}/k]),$$

where $i_1 j_1 \equiv 1 \pmod{k}$. The decomposition

$$\widetilde{L} = \mathbb{S}^{2n-1} = \mathbb{S}^1 * \cdots * \mathbb{S}^1 = \widetilde{L}(k; i_1) * \cdots * \widetilde{L}(k; i_n)$$

gives $\widetilde{L}$ a $CW$-structure, and

$$C(\widetilde{L}) \cong \mathbb{Z}[\mathbb{Z}/k] \xrightarrow{\cdot g^{j_n} - 1} \mathbb{Z}[\mathbb{Z}/k] \xrightarrow{\cdot \Sigma} \cdots \xrightarrow{\cdot g^{j_2} - 1} \mathbb{Z}[\mathbb{Z}/k] \xrightarrow{\cdot \Sigma} \mathbb{Z}[\mathbb{Z}/k] \xrightarrow{\cdot g^{j_1} - 1} \mathbb{Z}[\mathbb{Z}/k].$$

After tensoring with $\Lambda$ the complex is acyclic and $\cdot\Sigma$ is zero, so the other maps must be isomorphisms, and

$$\tau(C(L;\Lambda)) = \prod_{\ell=1}^{n} g^{j_\ell} - 1.$$

To complete the argument, some (not so sophisticated) algebraic number theory comes in. $\Lambda \cong \prod_{d|k} \mathbb{Q}[\zeta_d]$ where the product is over all divisor of $n$ greater than 1, and $\zeta_d$ is a primitive $d$-th root of 1. One then considers the quotient of the torsion of $L$ and the torsion of $L'$ as units in the corresponding number rings, and the Franz independence lemma (whose proof is similar to that of the Dirichlet Unit Theorem) says that the only way these units can be roots of unity is if the conditions of part (v) of Theorem 3.10 are satisfied. This gives the classification of classical lens spaces.

### 3.3.3. General remarks

The reason why we are spending so much time on spherical space forms is that this problem represents a paradigm for the construction and classification of finite group actions on manifolds. The existence discussion above moved from algebraic information to geometric information at the level of $CW$-complexes; this passage is accomplished via the Hurewicz theorem. The finiteness obstruction in $K_0$ measures when a finite $CW$-complex can be obtained. The torsion in $K_1$ allows the discussion of simple homotopy theory and allows the classification of classical lens spaces. It also plays a key role in transformation groups via the $s$-cobordism theorem.

What is missing? Two things – at least – the transition to non-free actions and the transition to manifold theory. The generalization of the homotopy theoretic techniques to the non-free equivariant case is in reasonable shape. For example, de Rham used torsion to classify linear but non-free, actions of finite groups on spheres. The general theory of the equivariant finiteness obstruction, equivariant simple homotopy theory, and the equivariant $s$-cobordism theorem is worked out (see [83]) at least for the $PL$ and smooth cases; these issues for locally linear actions on topological manifolds are still active areas of current research. The passage from algebraic information to geometric information at the level of $CW$-complexes in the non-free case was studied in [104] and [10].

The second missing ingredient in our discussion of the space form problem is the transition to manifolds. We have not yet addressed the existence question of when a Swan complex has the homotopy type of a manifold or the uniqueness question of classifying *all* manifolds within a homotopy type. Two aspects of these questions have been resolved nicely: the Madsen-Thomas-Wall result which gives the class of finite groups which arise as fundamental groups of manifolds covered by the sphere and Wall classification of fake lens spaces (i.e. manifolds homotopy equivalent to classical lens spaces) of dimension greater than five. Results on manifolds are usually accomplished via surgery theory and we discuss this next.

The systematic method for classification of manifolds is called *surgery*. The idea is that surgery theory reduces classification questions to a mix of algebraic topology and the algebra of quadratic forms. Some of the ingredients necessary for this reduction are handlebody theory, bundle theory, transversality, and embedding theorems. The embedding theorems make the theory most effective when the dimension of the manifold is $\geqslant 4$, where there is sufficient room to mimic algebraic manipulations by geometric embeddings. Both transversality and embedding theory provide "flies in the ointment" for the development of equivariant surgery. Hence, we concentrate on classical surgery theory.

As far as references for surgery theory, to the great detriment of the subject, there is no modern account. The most comprehensive is the book of Wall [136]. We refer the read to the paper of Milnor [95] and the book of Browder [25] for geometric background. For background on classifying spaces and bundle theory see Milgram-Madsen [92]. For information specific to spherical space forms, see Davis-Milgram [45], Wall [136] and Madsen-Thomas-Wall [86]. Modern aspects of surgery theory can be found in the books of Ranicki [113] and Weinberger [142], however they were not written with classical surgery theory as their main focus.

### 3.3.4. *Existence of space forms*

When does a finite group $G$ act freely on $\mathbb{S}^n$? As above, one must have $H^{n+1}(G, \mathbb{Z}) \cong \mathbb{Z}/|G|$. Then there exists a Swan complex, that is, a $CW$-complex $X$ with dim $X = n$, $\pi_1(X) = G$, and $\widetilde{X} \simeq \mathbb{S}^n$. In fact, the $k$-invariant gives a one-to-one correspondence between polarized homotopy types and additive generators of $H^{n+1}(G, \mathbb{Z})$. Any Swan complex is finitely dominated, and the finiteness obstruction $[X] \in \widetilde{K}_0(\mathbb{Z}G)$ vanishes if and only if $X$ has the homotopy type of a finite $CW$-complex. This can be effectively computed, [137], [91], [42], although the algebraic number theory can be quite involved. One qualitative result is that if $G$ is a group of period $n + 1$, there is always a finite Swan complex of dimension $2n + 1$.

To see when a finite Swan complex has the homotopy type of a manifold, one uses surgery theory, which provides necessary and sufficient conditions in dimensions greater than 4, and provides necessary conditions all dimensions. First note that $X$ is a *Poincaré complex*, i.e. $X$ is a finite complex and there exists a class $[X] \in H_n(X)$ so that $\cap \operatorname{tr}[X] : H^i(X; \mathbb{Z}G) \to H_{n-i}(X; \mathbb{Z}G)$ is an isomorphism for all $i$.

The first obstruction to a Poincaré complex having the homotopy type of a manifold is the existence of a lift of the Spivak normal bundle to $BTOP$. The Spivak normal bundle is the homotopy theoretic analogue of the stable normal bundle of a manifold and is defined as follows. Embed $X$ in $\mathbb{R}^K$ ($K$ large) and let $N(X)$ be a closed regular neighborhood. Convert the map $\partial N(X) \to X$ to a fibration $p : E \to X$, then it is a formal consequence of Poincaré duality that the fiber has the homotopy type of $\mathbb{S}^{K-n-1}$ (see [25]). Fibrations with fibers having the homotopy type of a sphere are called *spherical fibrations*. Spherical fibrations over $X$ are classified by $[X, BG]$ where $G = \operatorname{colimit} G_k$ is the stabilization of the topological monoid $G_k$ of self-homotopy equivalences of $\mathbb{S}^{k-1}$. The map $p : E \to X$ and its classifying map $X \to BG$ are both referred to as the *Spivak normal bundle of $X$*. The next step is to see whether the Spivak bundle lifts to an honest topological sphere bundle, classified by a map $X \to BTOP$, and, if so, in how many ways. The obstruction to lifting the Spivak bundle to $BTOP$ is an element in $[X, B(G/TOP)]$. For a general Poincaré complex this might be non-trivial, but for a finite Swan complex it vanishes. The argument is as follows. The space $B(G/TOP)$ is an infinite loop space (i.e. the 0-space of an $\Omega$-spectrum), and hence the abelian group $[X, B(G/TOP)]$ injects into $\oplus_p [X, B(G/TOP)]_{(p)}$ which injects into $\oplus_p [\widetilde{X}/G_p, B(G/TOP)]_{(p)}$, where $G_p$ is a $p$-sylow subgroup of $G$ and the second injection is via the transfer map. Now $\widetilde{X}/G_p$ has the homotopy type of a lens space or a quaternionic space form, so the Spivak bundle map to $B(G/TOP)$ vanishes.

If the Spivak bundle $\nu : X \to BG$ of a Poincaré complex $X$ lifts to $\tilde\nu : X \to BTOP$, then one can apply a transversality construction as a first step in the attempt to construct a manifold having the homotopy type of $X$.

**Definition 3.5.** A *Thom invariant* for $X$ is an element $\alpha \in \pi_{n+k}(T(\tilde\nu))$ for some lift $\tilde\nu : X \to BTOP(k)$ of the Spivak bundle $\nu$, so that $h(\alpha) \cap U = [X]$ where $h$ is the Hurewicz map and $U \in H^k(T(\tilde\nu))$ is the Thom class.

Given any lift $\tilde\nu$, Thom invariants always exist (see [25]) and are essentially given by collapsing out the complement of a regular neighborhood of $X$ in $\mathbb{S}^{n+k}$. Given a Thom invariant $\alpha : \mathbb{S}^{n+k} \to T(\tilde\nu)$, one may take the complement of the 0-section $X$ to obtain a degree one map

$$f : M = \alpha^{-1}(X) \to X$$

where $M$ is a closed manifold. Furthermore, transversality gives a trivialization of $f^*(\tilde\nu) \oplus \tau_M$. Hence we call the induced map

$$(f, \hat{f}) : (M, \nu_M) \to (X, \tilde\nu)$$

a *degree one normal map*, where $\nu_M$ is the stable normal bundle of $M$ (equipped with a trivialization of $\nu_M \oplus \tau_M$). For such a map, if $[\beta] \in \ker(\pi_i M \to \pi_i X)$ is represented by an embedding $\beta : \mathbb{S}^i \to X$ one may use the normal data $\hat{f}$ to thicken $\beta$ up to an embedding $\mathbb{S}^i \times D^{n-i} \to M$ and perform surgery to obtain

$$M' = (M - \mathbb{S}^i \times \text{int } D^{n-i}) \cup_{\mathbb{S}^i \times \mathbb{S}^{n-i-1}} (D^{i+1} \times \mathbb{S}^{n-i-1})$$

and a degree one normal map $(g, \hat{g}) : (M', \nu_{M'}) \to (X, \tilde\nu)$. (In effect you are killing the homotopy class $[\beta]$, see [95].) The equivalence relation generated by surgery is called normal bordism and the set of equivalence classes is denoted $N(X)$. The main theorem of surgery theory is that there is the surgery obstruction map

$$\theta : N(X) \to L_n(\mathbb{Z}\pi_1 X).$$

$L_n(\mathbb{Z}\pi_1 X)$ is a algebraically defined abelian group associated to the group ring, closely related to quadratic forms over the group ring. If $(f, \hat{f})$ is normally bordant to a homotopy equivalence, then $\theta(f, \hat{f}) = 0$, and conversely, for $n > 4$, if $\theta(f, \hat{f}) = 0$, then $(f, \hat{f})$ is normally bordant to a homotopy equivalence.

An application to the spherical space form problem is given by:

**Theorem 3.11** (Petrie)**.** *Let $G = \mathbb{Z}_7 \times_T \mathbb{Z}_3$ be the semidirect product where $\mathbb{Z}_3$ acts non-trivially on $\mathbb{Z}_7$. Then $G$ acts freely on $\mathbb{S}^5$, but does not act freely and linearly.*

**Proof:** It is an easy computation to show that $G$ can't act freely and linearly on $\mathbb{S}^5$. The Lyndon-Hochschild-Serre spectral sequence shows that $H^6(G, \mathbb{Z}) \cong \mathbb{Z}/21$, so $G$ has period 6 and there is a Swan complex of dimension 5. Since $T(G) = \tilde{K}_0(\mathbb{Z}G)$ (see [137, p. 545-546]), there is a Swan complex with a zero finiteness obstruction, hence there is a finite Swan complex $X$ of dimension 5 and fundamental group $G$. We have already observed that the Spivak bundle of $X$ lifts to $BTOP$, so there is a degree one normal map

$$(f, \hat{f}) : (M, \nu_M) \to (X, \xi).$$

Now since $G$ is odd order and $n$ is odd, $L_n(\mathbb{Z}G) = 0$ (see [16]), one may complete surgery to a homotopy equivalence $M' \to X$. By the generalized Poincaré conjecture, $M' \cong \mathbb{S}^5/G$ for a free $G$-action of $M'$. (Note: a similar analysis gives a free smooth action on $\mathbb{S}^5$.)

**Remark 3.1.** *Petrie [107] proved his theorem in a much more explicit and elementary manner. He noted that there is a free $G$-action on the Brieskorn variety*

$$\Sigma^5 = \{(z_1, z_2, z_3, z_4) \in \mathbb{C}^4 \mid z_1^7 + z_2^7 + z_3^7 + z_4^3 = 0\} \cap \mathbb{S}^7.$$

*This is only a rational homology sphere, but Petrie shows how to perform surgery to make it into a sphere. Another application of this sort of technique is given in [44].*

The analysis of the problem of determining when a group $G$ of period $n > 5$ can act freely on a sphere is much more difficult. One must compute surgery obstructions for all $k$-invariants and normal invariants. We state the beautiful result of Madsen-Thomas-Wall.

**Definition 3.6.** Let $n$ be a positive integer. A group $G$ satisfies the *$n$-condition* if every subgroup of order $n$ is cyclic.

**Theorem 3.12** (Madsen-Thomas-Wall [86])**.** *A finite group $G$ acts freely on some sphere if and only if $G$ satisfies the $p^2$- and $2p$-conditions for all primes $p$.*

In contrast, Wolf [143] analyzed the free, linear actions on spheres:

**Theorem 3.13.** *A finite group $G$ acts freely and linearly on some sphere if and only if the following two statements are satisfied.*
  (i) *$G$ satisfies the $pq$-condition for all primes $p$ and $q$.*
  (ii) *$G$ has no subgroup isomorphic to $SL_2(\mathbb{F}_p)$ for a prime $p > 5$.*

*3.3.5. Uniqueness of space forms*

Surgery theory also attacks the uniqueness question – classifying manifolds within a homotopy type. It was motivated by the proof of the generalized Poincaré conjecture and the Kervaire-Milnor classification of exotic spheres. Two great successes of surgery theory are various cases of the Borel conjecture (especially that there are no fake tori of dimension greater than 3), and the classification of fake lens spaces due to Browder, Petrie, and Wall.

**Definition 3.7.** A *fake lens space* is a space with covered by the sphere with cyclic fundamental group.

Every fake lens space is homotopy equivalent to a lens space.

**Theorem 3.14** ([136])**.** *Two fake lens spaces with odd order fundamental group and dimension greater than 3 are homeomorphic if and only if they have the same Reidemeister torsion and $\rho$-invariants.*

Furthermore, Wall shows exactly which invariants are realized. Reidemeister torsion was defined in Section 3.3.2. The $\rho$-invariant is defined as follows.

**Definition 3.8.** For a closed, oriented, odd-dimensional manifold $M$ with finite fundamental group $G$, define the *$\rho$-invariant*

$$\rho(M) = \frac{1}{s}\mathrm{sign}(G, W) \in \widetilde{K}_0(\mathbb{C}G) \otimes \mathbb{Q},$$

where $s$ disjoint copies of the universal cover $\widetilde{M}$ bound a compact, oriented, free $G$-manifold $W$.

Given an action of a finite group $G$ on a compact, oriented, even-dimensional manifold $X$, Atiyah-Singer [15, p. 273-274], defined the $G$-signature $\mathrm{sign}(G, X) \in \widetilde{K}_0(\mathbb{C}G)$. We resist the urge to state the definition and to indicate the many applications in transformation groups of the Atiyah-Singer G-signature theorem.

*3.4. Final Remarks*

In this paper we have attempted to cover the main ideas and examples which have made the subject of transformation groups an important and highly developed subject in topology. The background required to work on any remaining questions may require a combination of skills in manifold theory, group cohomology, algebraic $K$- and $L$-theory, homotopy theory and differential geometry. However daunting this may be, it is apparent that much remains to be understood about topological symmetries and we hope that the reader of this paper will take it upon himself to explore this topic further.